\documentclass[11pt,reqno,oneside]{amsart}
\usepackage{color}
\usepackage{amssymb,amsmath,mathrsfs,amsthm}
\usepackage{caption}
\usepackage{enumitem}
\usepackage{mathtools}
\usepackage{geometry}
\usepackage{amsfonts}
\usepackage{color}

\usepackage{esint,comment}
\usepackage{apptools}
\AtAppendix{\counterwithin{thm}{section}}

\usepackage{graphicx}
\allowdisplaybreaks
\newtheorem{thm}{Theorem}
\newtheorem{cor}[thm]{Corollary}

\newtheorem{prop}[thm]{Proposition}

\theoremstyle{definition}

\newtheorem{rem}{Remark}

\def \no#1#2#3 {{\bf #1} (#3), #2.}
\def \eds#1#2#3 {#1, #2, #3.}

\def\R{{\mathbb R}}

\def\d{{\rm d}}

\def\:{{\colon}}

\def\be#1{\begin{equation}\label{#1}}
\def\ee{\end{equation}}

\def\<{\langle}
\def\>{\rangle}
\def\coloneqq{:=}

\newcommand{\lec}{\lesssim}
\newcommand{\bs}{\begin{split}}
\newcommand{\essss}{\end{split}}

\newcommand{\eqnb}{\begin{equation}}
\newcommand{\eqne}{\end{equation}}

\newcommand{\es}{\mathcal{S}}
\newcommand{\uu}{\mathcal{V}}

\renewcommand{\ee}{\mathrm{e}}

\newcommand{\p}{\partial}

\newcommand{\re}{\mathrm{Re}}
\newcommand{\im}{\mathrm{Im}}

\renewcommand{\R}{\mathbb{R}}

\newcommand{\C}{\mathbb{C}}

\renewcommand{\d}{\mathrm{d}}

\begin{document}

\title[Delta wings of small aspect ratio]{An analysis of aerodynamic properties of delta wings} 
\author{W. S. O\.za\'nski}
\address{Florida State University, Tallahassee, FL 32306, and Princeton University, Princeton, NJ 08540}
\email[]{wozanski@fsu.edu}
\maketitle

\date{}

\medskip

\begin{abstract}
We consider a sharp-edge delta wing of small aspect ratio $A>0$, which is an example of a 3D airfoil whose aerodynamic properties cannot be modeled using the \emph{potential lift} only. An important role is played by a pair of attached vortices, which generate the \emph{vortex lift}. We review in detail the leading-edge suction analogy, developed by Edward Polhamus (\emph{NASA Technical Note D-3767}, 1966), which shows that the  most accurate description of the vortex lift can be obtained by avoiding \emph{any} direct analysis of  vorticity of the flow around the wing. We also describe a 2D slender wing theory, which was introduced by Jones (\emph{NACA report no.~835}, 1946) to approximate the potential lift. We show that, neglecting certain two borderline nonintegrability issues and an elliptic regularity issue, the slender wing theory is sufficient to yield accurate prediction of the potential lift as $A\to 0$. We also show that this 2D approximation is necessary to justify the modelling assumptions of the leading-edge suction analogy, so that the slender wing model can also be used to approximate the vortex lift as $A\to 0$. For completeness, we provide additional explanations of the necessary airfoil theory.
\end{abstract}



\section{Introduction}\label{sec_intro}

\begin{figure}[htbp]
\centering
    \includegraphics[width=7cm]{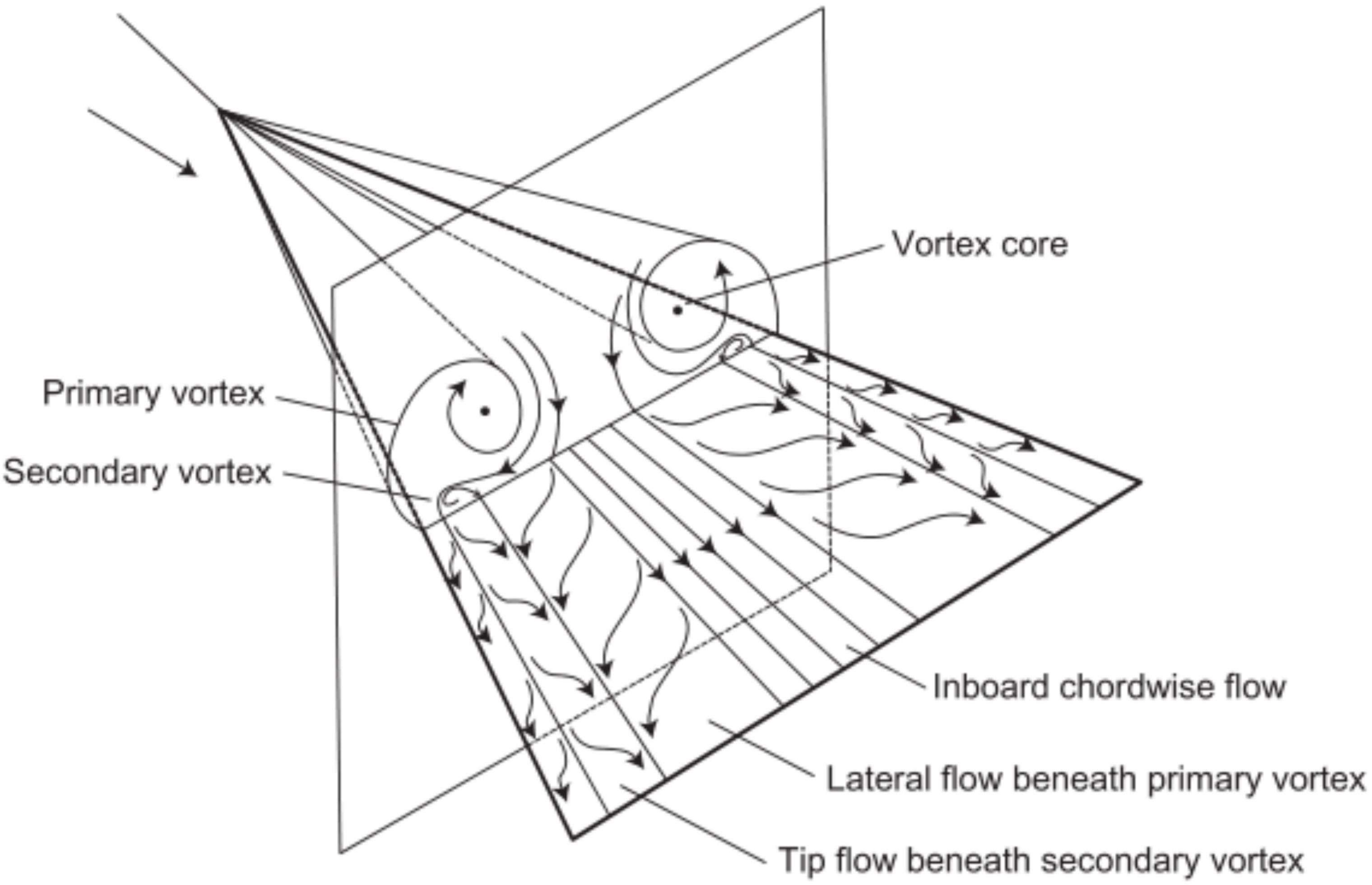}
    \hspace{0.1cm}
      \includegraphics[width=4cm]{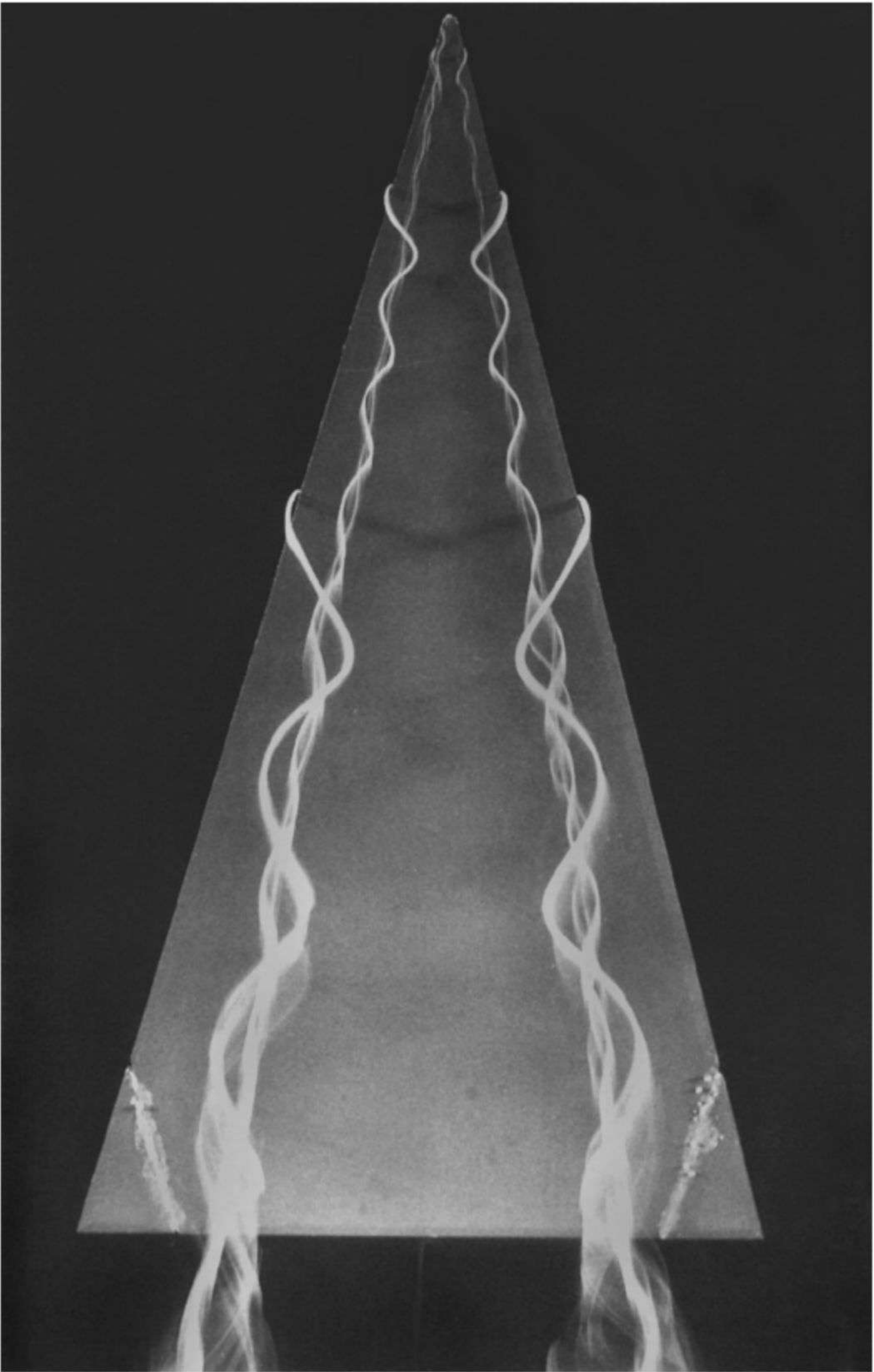} \hspace{0.1cm}
      \includegraphics[width=3.5cm]{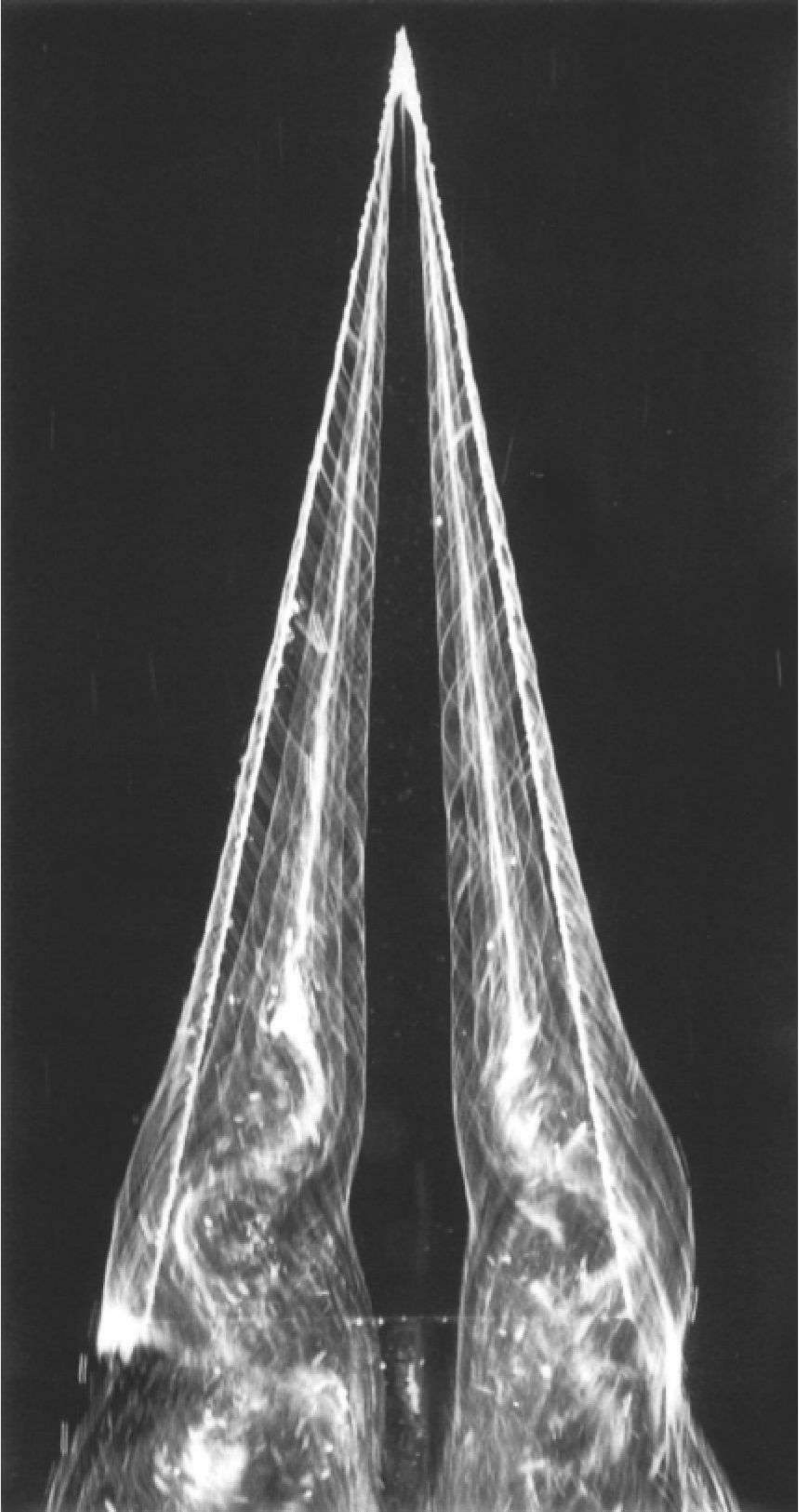}
  \caption{Left (after~\cite[Fig.~5.33]{houghton}): A sketch of the attached vortices of sharp-edge delta wings. Middle (after~\cite[Fig.~8]{werle} and \cite[Fig.~90]{van_dyke}): A snapshot of the vortices (visible by a colored fluid in water) above a sharp-edge delta wing with $\beta = 15^\circ$ and angle of attack $\alpha = 20^\circ$. Right (after \cite{solignac}): A snapshot of attached vortices with vortex breakdown (occurring at about 0.7 maximum chord), visible by electrolysis in water flow. Here $\alpha = 35^\circ$ and Reynolds number is $3000$. }\label{fig_sketch_vortices}
\end{figure}

We consider a 3D airfoil of the form a sharp-edge delta wing $\es$ of a fixed area $|\es|$. We denote  the maximal wingspan of the wing by $b>0$ and the aspect ratio by
\eqnb
A \coloneqq\frac{b^2}{|\es|}= 4 \tan \beta ,
\eqne
where $2 \beta \in (0,\pi)$ denotes the angle at the leading tip of the wing, see Fig.~\ref{fig_sketch}. Such 3D airfoils cannot be treated by extending the standard 2D airfoil theory, particularly for small aspect ratio $A$.  It is well-known that such wings generate attached vortices, see Fig.~\ref{fig_sketch_vortices}. These are in fact vortex columns, as the axial velocity varies with the radius (as can be observed from Fig.~\ref{fig_sketch_vortices}(right)), and are interesting from the perspective of hydrodynamic instability. For example, Gallay \& Smets \cite{gs_linear,gs_spectral} proved linear and spectral instability of certain vortex columns under axisymmetric perturbations, while \cite{AO} provieded rigorous construction of infinitely many modes of instabilities of certain vortex columns (which include Batchelor trailing vortices) of the form of ``ring modes'', see also \cite{ls,cs,ls_87,hg}.  The attached vortices also exhibit vortex breakdown, see Fig.~\ref{fig_sketch_vortices}.  We refer the reader to an enlightening exposition of aeronautical developments related to delta wings by Luckring~\cite{luckring}.

\begin{figure}[htbp]
\centering
    \includegraphics[width=11cm]{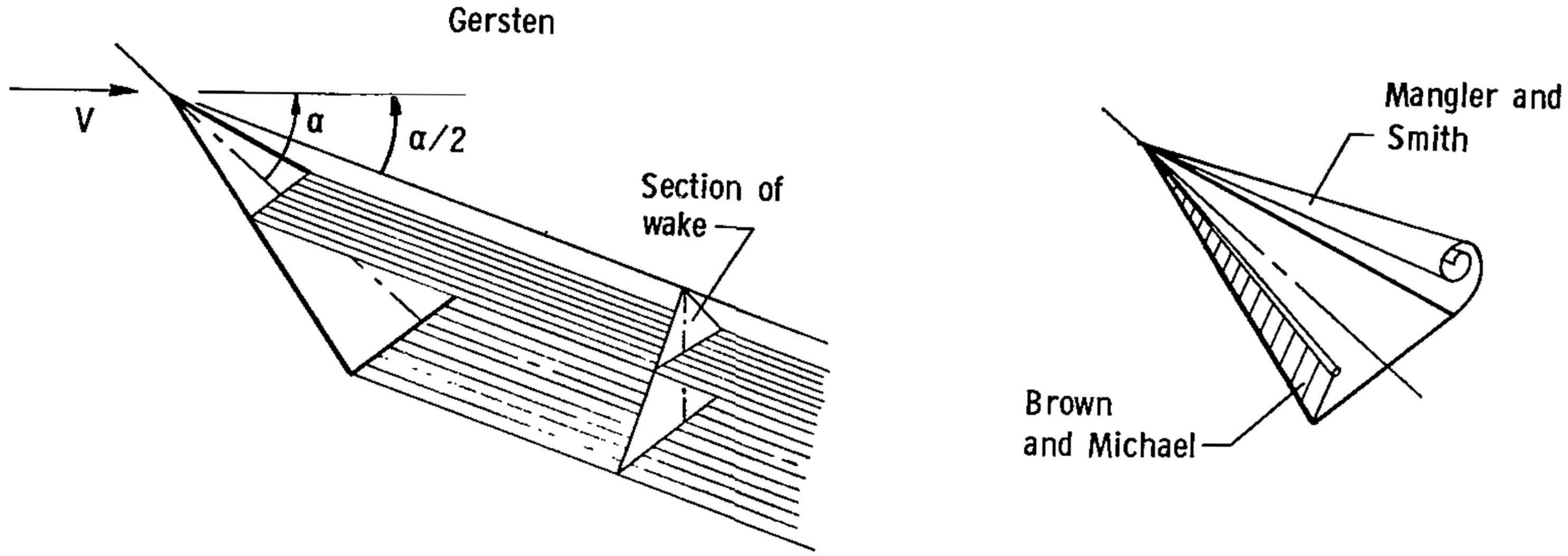}       
  \caption{(after \cite[Fig.~4]{polhamus}) A sketch of the approaches of \cite{brown_michael,mangler_smith,gersten} (see also \cite{legendre_52}) to describe the vortex lift $L_v$.}\label{fig_approaches}
\end{figure}

 Interestingly, the appearance of the attached vortices also affects aerodynamic properties of delta wings. Indeed, apart from the \emph{potential lift} $L_p$, which arises from a potential flow around the wing and its wake (which can be computed using a solution to the lifting surface problem, see Section~\ref{sec_LSP} for details), additional lift, the \emph{votex lift} $L_v$, is generated by the vortices. At first sight it appears that $L_v$ can be computed by adapting a 2D model of point vortex suction (which we describe in detail in Appendix~\ref{app_2d_vortex_suc}). This approach was applied by Legendre~\cite{legendre_52}, but the theoretical approach turned out to be in poor agreement with experiment (see~Fig.~7 in \cite{gersten}, for example). Later Brown and Michael \cite{brown_michael} considered a model consisting of a vortex filament and a feeding flat vortex sheet, Mangler and Smith \cite{mangler_smith} developed a model of a spiraling vortex sheet, and  Gersten \cite{gersten} considered a three-dimensional vortex wake approach, see Fig.~\ref{fig_approaches} for a sketch. Remarkably, these attempts provide very limited agreement with experiments (see Fig.~\ref{fig_results}(right)).

\begin{figure}[htbp]
\centering
    \includegraphics[width=5cm]{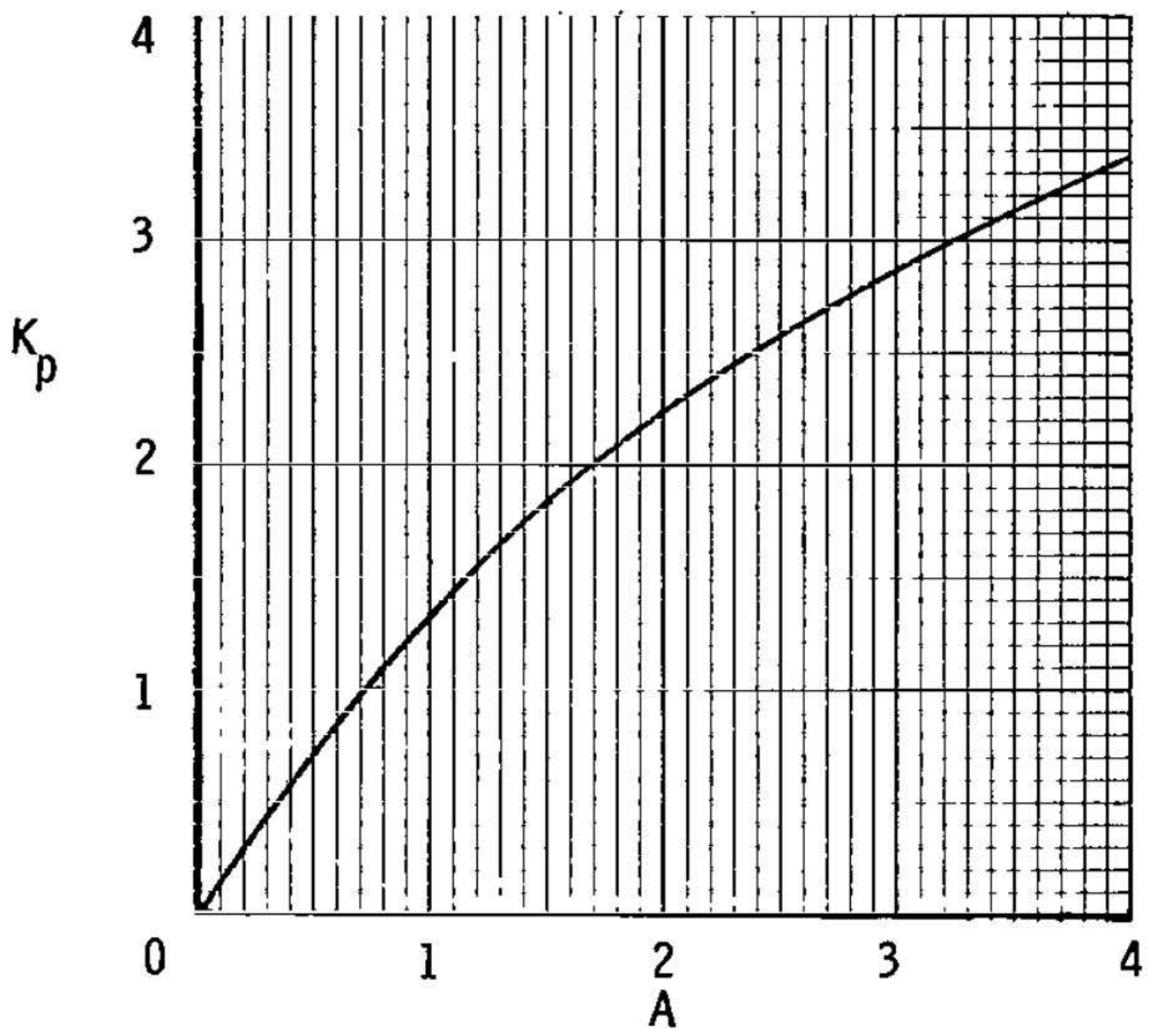}
    \hspace{1cm} 
    \includegraphics[width=5cm]{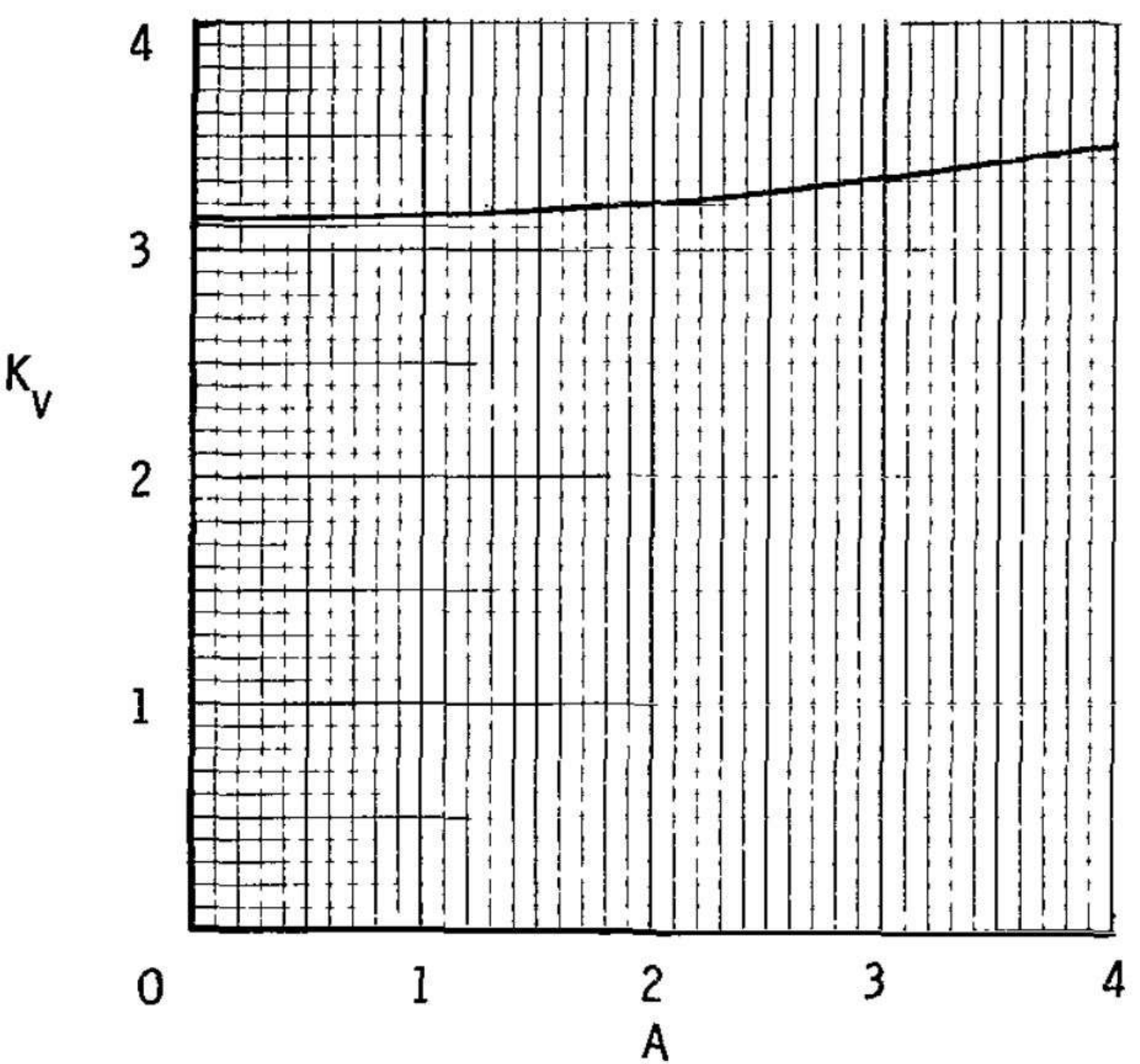}    
  \caption{(after \cite[Fig.~6,9]{polhamus}) Numerical values of aerodynamic constants $K_p$, $K_v$ computed by Polhamus \cite{polhamus}.}\label{fig_KpKv}
\end{figure}

The most effective way of describing it is the edge-suction analogy introduced by Polhamus~\cite{polhamus}, whose model argued that total lift consists of the \emph{potential lift} $L_p$ and the \emph{vortex lift} $L_v$,
\eqnb\label{polhamus_ansatz}
\begin{split}
L & = { L_p} + {L_v}\coloneqq  \frac12 \rho V^2 |\es| \left( { K_p \sin \alpha \cos^2 \alpha}  + {K_v \cos \alpha \sin^2 \alpha } \right),
\end{split}
\eqne
where $K_p $, $K_v$ are a dimensionless constants, depending on the aspect ratio $A$ only, see Fig.~\ref{fig_KpKv}. Also, \cite{polhamus} represents $K_v$ as
\eqnb\label{Kv}
K_v \coloneqq \left( K_p - K_p^2 K_i \right) \frac{1}{\cos \Lambda},
\eqne
where $K_i \coloneqq \frac{\p C_{D_i}}{\p C_L^2 }$ is the induced drag to lift ratio (see Section~\ref{sec_polhamus} details). 
The lift formula \eqref{polhamus_ansatz} is based on the central modelling assumption of the \emph{leading-edge-suction analogy}, which we explain in detail in Section~\ref{sec_polhamus_explanation}. Moreover, as mentioned in \cite[p.~7]{polhamus} the values of $K_p$ and $K_v$ (as shown in Fig.~\ref{fig_KpKv}) are ``determined by a modification of the  Multhopp \cite{multhopp}  lifting-surface theory'', which describes a numerical algorithm for solving \eqref{LSP}. The leading edge-suction analogy turns out to be in a remarkable agreement with experiments, see Fig.~\ref{fig_results}, and was very successful in developing variable-sweep aircrafts\footnote{Edward Polhamus had a very rich career at NACA/NASA Langley research center, and the leading-edge-suction analogy became a central idea leading to the development of variable sweep military aircraft  in the second half of the twentieth century, see \cite{polhamus_84} as well as the 1945 paper~\cite{jones_45} of Jones, who showed that the delta wings with angles smaller than the Mach cone are most suitable for supersonic flights. See also Polhamus's enlightening overview~\cite{polhamus_86} of the vortex lift research, as well as Luckring's article~\cite{luckring_2016} on scientific contributions of Polhamus.}.

The lifting surface problem (which we describe in detail in Section~\ref{sec_LSP}) is concerned with finding a potential function $\Phi$ describing the potential velocity field $\nabla \Phi$ around the wing and its wake. It is a major open problem of mathematical aerodynamics; so far only numerical algorithms, such the Multhopp~\cite{multhopp} algorithm mentioned above, and the refinements   due to Van Spiegel and Wouters~\cite{van_spiegel_wouters} and Lamar~\cite{lamar},  provide any understanding of the problem.

\begin{figure}[htbp]
\centering
    \includegraphics[width=12cm]{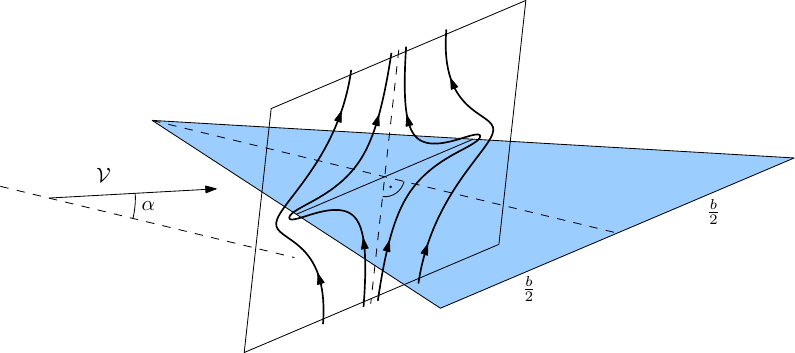}
  \caption{A sketch of the 2D potential flow  approximation of the flow around a delta wing.}   \label{fig_2d}
\end{figure}

\begin{figure}[htbp]
\centering
    \includegraphics[width=10cm]{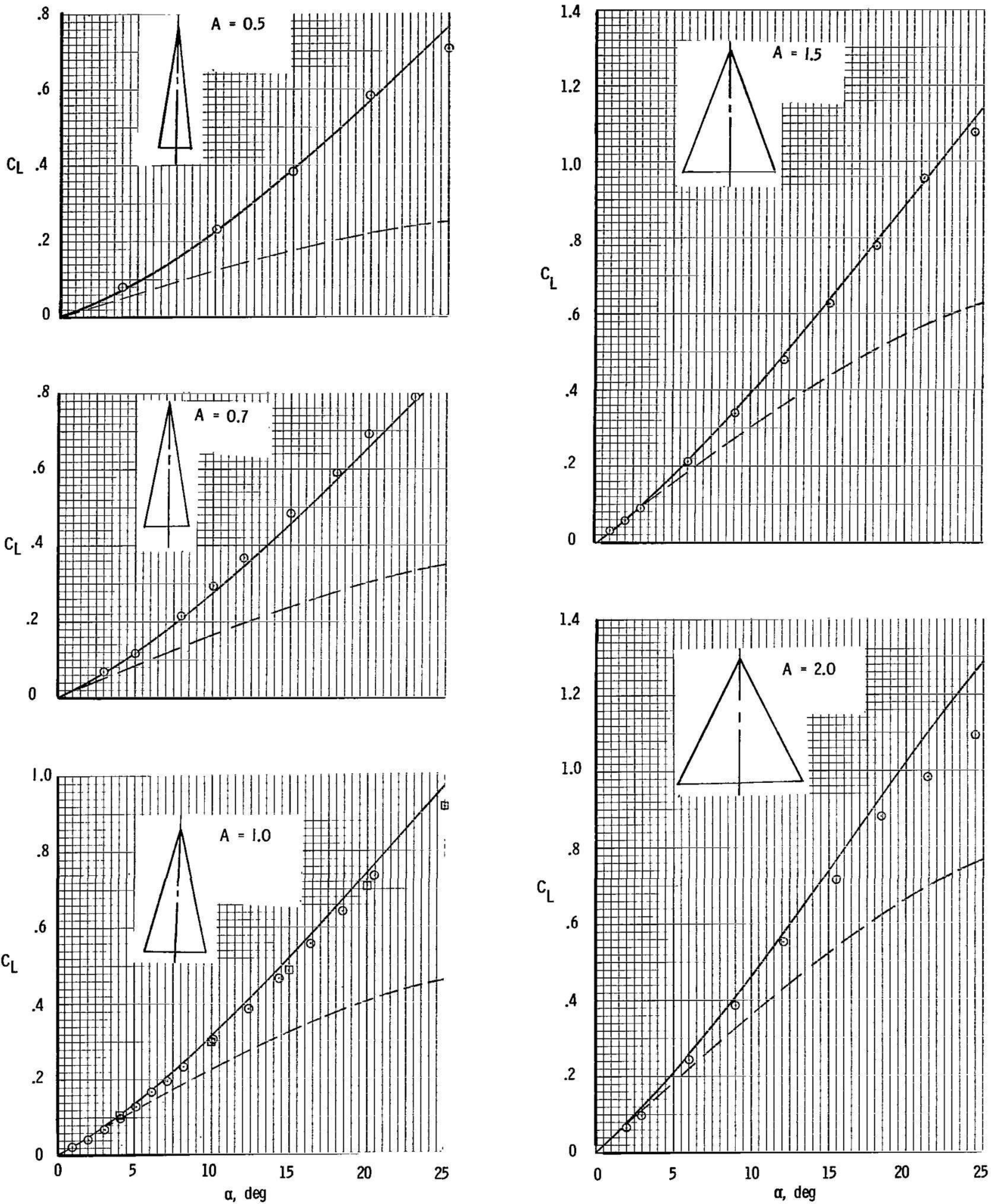}   \hspace{0.1cm} 
    \includegraphics[width=6cm]{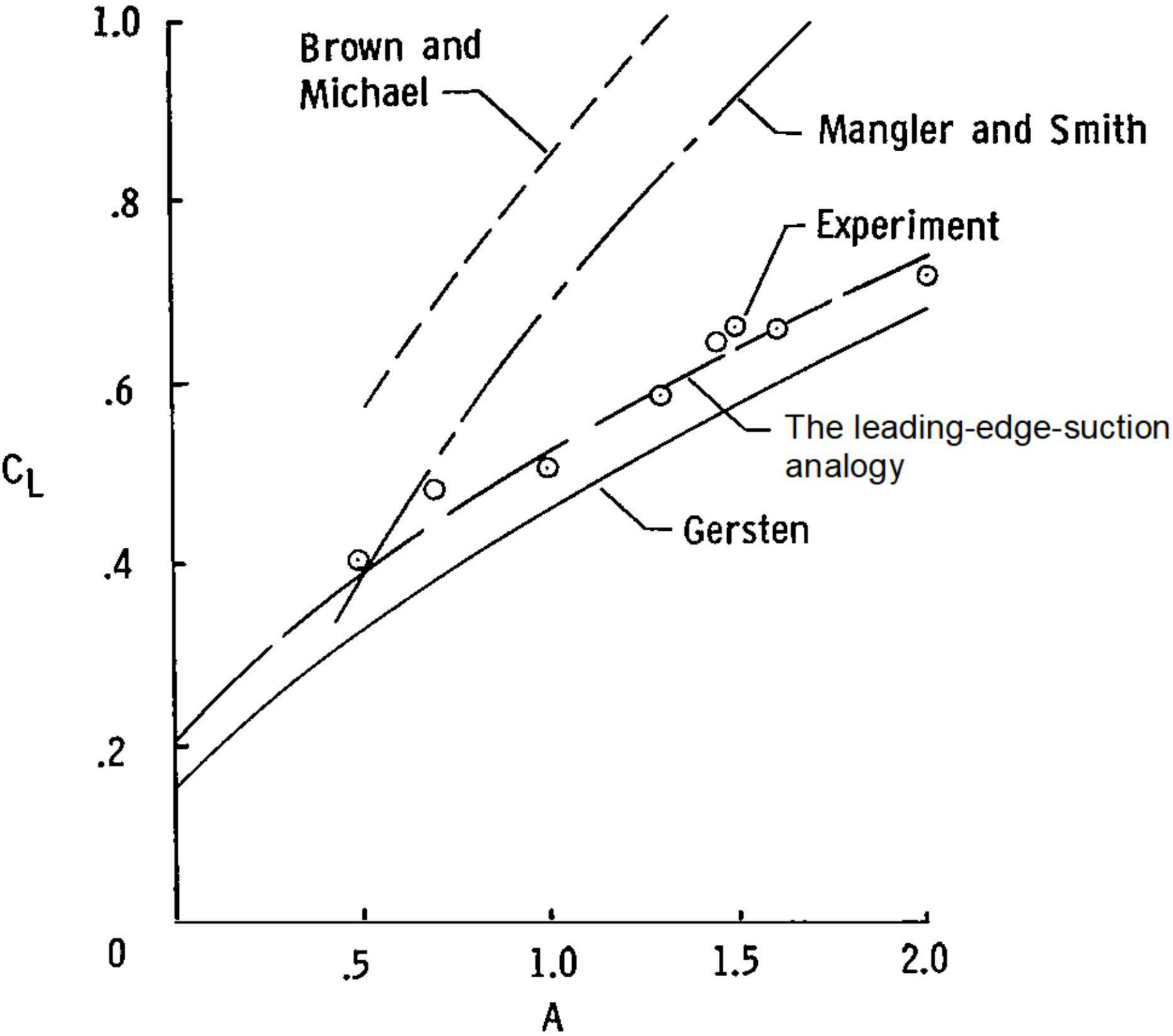} 
  \caption{Results of Polhamus's approach. Left (after~\cite[Fig.~12]{polhamus}):   Comparison of the leading edge-suction analogy \eqref{polhamus_ansatz} with experiment. Here the dashed line corresponds to $C_{L_p} \coloneqq \frac{L_p }{\frac12 \rho V^2 |\es|}$ only, the solid line to $C_L\coloneqq \frac{L_p + L_v }{\frac12 \rho V^2 |\es|}$, and the experimental data  is from \cite{ref4,ref11,ref12,ref13}. Right  (after~\cite[Fig.~11]{polhamus}): Comparison  of \eqref{polhamus_ansatz} to the vorticity-based approaches of \cite{legendre_52,brown_michael, gersten,mangler_smith}. Here   $\alpha = 15^\circ$.}\label{fig_results}
\end{figure}

A solution $\Phi$ to the lifting surface problem \eqref{LSP} determines the potential part $L_p$ of the lift by integrating the jump of $\p_1 \Phi$ across the wing surface (see Section~\ref{sec_LSP} for details), and so determines, to the leading order, the basic aerodynamic properties of \emph{any flat wing}. For the particular case of the delta wings of small aspect ratio $A$, some modern expositions replace $\Phi$ by the potential $\phi$ of a 2D flow in a plane perpendicular to the wing;  namely, the fluid flow is approximated by the potential flow $\nabla \phi$ around a flat line, see \eqref{phi_def} and Appendix~\ref{app_2dpot} for details, and Fig.~\ref{fig_2d} for a sketch. This is an interesting approach going back to the work of Jones~\cite{jones}, but the approximation  lacks any rigorous justification. This is not surprising, since the lifting surface problem is unresolved  in the case of the delta wing. 

The purpose of this note is twofold. The first purpose is to better understand the formula \eqref{polhamus_ansatz}. To this end, we clarify (in Section~\ref{sec_polhamus_explanation}) what \emph{exactly} is meant by Polhamus~\cite{polhamus} in the derivation of \eqref{polhamus_ansatz}, particularly regarding the algorithm for finding $K_p$ in Fig.~\ref{fig_KpKv}, the leading-edge suction analogy (between the suction force $S$ and the vortex lift $L_v$), and the derivation of $K_v$ as in \eqref{Kv}. Furthermore, the treatment of $K_i$, which is left implicit in the original work \cite{polhamus}, and we discuss its key role in understanding the dependence of $K_v$ on the aspect ratio $A$ as plotted  in Fig.~\ref{fig_KpKv}. It is particularly nuanced, as it relies on a \emph{change of modelling assumptions}, the implicit assumption that the delta wings of small aspect ratio are optimal in the sense of minimizing the induced drag $D_i$ and that the downwash velocity field is constant throughout $\es$. We discuss these issues in Section~\ref{sec_polhamus_explanation} and  we show that the last two  issues \emph{can only be justified} by using the 2D potential flow mentioned above.

This brings us to the second purpose of this note, which is to recall and justify the 2D potential flow approach for delta wings of small aspect ratio. It turns out that analyzing the error between the two potentials, $\phi - \Phi$, naturally leads to three analytical issues: 
\begin{align}
&|x'|^{-1} \text{ does not belong to } L^2_{\rm loc } (\R^2) , \label{issue1}\\
&\text{the trace operator is not bounded from  } H^{1/2}(\R^3) \text{ onto }L^2 (\R^2), \label{issue2}\\
&\text{elliptic regularity \eqref{the_elliptic_issue} (of a Poisson problem in }\R^3\setminus (\es \cup \widetilde{\es})\text{) is missing}. \label{issue3}
\end{align}
Here $\widetilde{\es}$ denotes the wake \eqref{wake_def} of the wing $\es$.
 We  show that, if \eqref{issue1}--\eqref{issue3} are neglected, then the error can be controlled.
\begin{prop}[2D potential flow approximation]\label{P1}
If $\phi$ denotes the 2D solution (at each $x_1$) of the potential flow around the delta wing given by \eqref{2d_pot}, $\Phi$ denotes a solution of the ($3$D) lifting surface problem \eqref{LSP} in the case of the delta wing of aspect ratio $A>0$, then, neglecting the issues \eqref{issue1}--\eqref{issue3},
\eqnb\label{claim_prop}
\| \phi - \Phi \|_{H^{3/2} (\R^3 \setminus (\es \cup \widetilde{\es}) )} \lec V \sin \alpha \, A^{\frac32}  .
\eqne
\end{prop}
We note that an error estimate of this form is challenging, since it is not clear that any of $\Phi, \phi$ belong to $H^{3/2} (\R^3 \setminus (\es \cup \widetilde{\es}) )$, see~\eqref{LSP}, \eqref{regularity_of_phi} for the regularity of $\Phi$, $\phi$.  We note that  the $H^{3/2}$ topology of the error estimate in \eqref{claim_prop} is needed in order to approximate $K_p$, $K_v$ by the values resulting from the 2D approximation. Actually, we would really like to control $\phi -\Phi$ in $H^{3/2+\delta }$ (for any $\delta>0$), since, recalling that the potential lift is obtained by integrating $\p_1 \Phi$ over $\es$ (see \eqref{potential_lift} for details),  we could then estimate the  lift error using the trace operator,  $| \int_{\es} (\p_1 \phi - \p_1 \Phi ) | \lec_\delta \| \phi - \Phi \|_{H^{\frac32+\delta}}$. This is another example of the relevance of issue \eqref{issue2}. Nevertheless,  Proposition~\ref{P1} helps understand Fig.~\ref{fig_KpKv} through the following.

\begin{cor}[$K_p$ and $K_v$ asymptotics]\label{C1}
Suppose that the claim \eqref{claim_prop} is valid, neglect issue~\eqref{issue2} and accept Polhamus modelling assumption~\eqref{pol_change}. Then
\eqnb\label{Kp_and_Kv_claim}
K_p = \frac{\pi }2 A + O( A^{3/2} ), \qquad K_v = \pi + O(A^{1/2})
\eqne
as $A\to 0$.
\end{cor}
We note that the asymptotic expansions \eqref{Kp_and_Kv_claim} (which are in agreement with Fig.~\ref{fig_KpKv}), are not mentioned in \cite{polhamus}, which only points out (below \cite[(13)]{polhamus}) that ``the value of $K_v$ increases slightly from a value of $3.14$ for an aspect ratio of $0$ to a value of about $3.45$ for an aspect ratio of $4$''\footnote{However, it is noted that $K_v=\pi$ at $A=0$ in Fig.~4 in a subsequent paper \cite{polhamus_71} of Polhamus, with no justification. It is also claimed in another work \cite[(6)]{polhamus_69} of Polhamus, that, in supersonic flows, $K_p=\frac{\pi A}{2E}$, where $E$ is a certain elliptic integral of second kind. This is claimed to follow from a 1946 derivation by Stewart~\cite{stewart_46}, see also \cite{jones_47}.}. 

A resolution of the issues \eqref{issue1}--\eqref{issue3} remains a very interesting open problem. Since the asymptotic expansions \eqref{Kp_and_Kv_claim} give errors of strictly higher order than the leading terms, we conjecture that a more careful analysis can provide unconditional rigorous results, see Remarks~\ref{rem_elliptic} and \ref{rem_borderline} for other comments. 

We introduce the lifting-surface problem in Section~\ref{sec_LSP}, which is followed by an analysis of the delta wings and a detailed derivation of \eqref{polhamus_ansatz} in Section~\ref{sec_polhamus}. In Section~\ref{sec_2d_approach} we discuss the 2D slender wing theory and we prove Proposition~\ref{P1} in Section~\ref{sec_pf_P1} and Corollary~\ref{C1} in Section~\ref{sec_pf_cor}. 

Finally, for completeness of the exposition, we also provide a number of clarifying remarks and figures throughout this note, as well as Appendices~\ref{app_2d_vortex_suc}--\ref{app_munk} describing in detail some well-established facts. 

\section{The lifting surface problem}\label{sec_LSP}
Let $\rho >0$ denote the fluid density. Let us consider the background flow $V$ at angle of attack $\alpha \in (0,\pi/2 )$, that is 
\eqnb\label{uu_def}
\uu \coloneqq (V \cos \alpha , 0 , V \sin \alpha ),
\eqne
and let us consider a velocity field $u$ such that $\uu + u$ satisfies the steady 3D incompressible  Euler equations,
\eqnb\label{3d_euler}
(\uu + u )\cdot \nabla (\uu + u ) = -\frac{1}\rho \nabla P.
\eqne
Now consider a closed flat surface $\es \subset  \{ x_3 =0 \} \subset \R^3$, and denote by 
\eqnb\label{wake_def}
\widetilde{\es} \coloneqq \{ (x_1,x_2, 0 ) \colon \exists_{z_1 < x_1 } (z_1,x_2,0) \in \es \} 
\eqne
its wake (see Fig.~\ref{fig_sketch} for the case when $\es$ is a delta surface).

We want to find $u $ describing the perturbation of the background flow $\uu$ resulting from placing $\es$ in the flow. We will only consider $u$ which is potential in $\R^3 \setminus (\es \cup \widetilde{\es })$, say 
\[
u = \nabla \Phi,
\]
in which case we have $u\cdot \nabla u = \p_i\Phi  \p_i \nabla \Phi = \nabla | \nabla \Phi |^2 /2$, so that the Euler equations \eqref{3d_euler} are equivalent to their linearization around $\uu$, 
\eqnb\label{3d_euler_prelin}
V\cos \alpha\, \p_1 u + V \sin \alpha \, \p_3 u = -\frac1\rho \nabla P
\eqne
in $\R^3 \setminus (\es \cup \widetilde{S})$. We also expect that $| u| = O(\alpha )$ as $\alpha \to 0$, and so we simplify \eqref{3d_euler_prelin} by neglecting the second term on the left-hand side to obtain 
\eqnb\label{lin_euler}
V \cos \alpha \p_1  u = - \frac{1}\rho \nabla P 
\eqne
in $\R^3 \setminus (\es \cup \widetilde{S})$. This shows that we should expect 
\eqnb\label{P_def}
P = - \rho V \cos \alpha \,\p_1 \Phi
\eqne
in $\R^3 \setminus (\es \cup \widetilde{S})$. Therefore, by the Kutta condition (that $\widetilde{\es }$ does not generate lift) we expect that $\p_1 \Phi$ should be continuous on $\R^3 \setminus \es$. Moreover, we also demand that $\p_3 \Phi $ is continuous across $\widetilde{\es}$, so that $u$ remains weakly divergence-free. Namely, 
\[0=\int_{\R^3 \setminus \es } u \cdot \nabla \varphi = -\int_{\widetilde{\es}} \left( \p_3 \Phi (x_1,x_2,0^+) - \p_3 \Phi (x_1,x_2, 0^- )\right) \varphi (x_1,x_2) \,\d (x_1, x_2)\]
for every $\varphi \in C^1 (\R^3 \setminus \es)$, where we integrated by parts in the second step. 

In summary, given the penetration velocity field $w\colon \es \to \R$ we aim to find  
 $\Phi  \colon \R^3 \to \R$ (the velocity potential) such that  
\eqnb\label{LSP}
\begin{cases}
&\Delta \Phi = 0  \qquad \text{ in } \R^3\setminus ( \es \cup \widetilde{\es} ),  \\
&\p_3 \Phi = w \qquad \text{ on }\es,\\
& \p_1 \Phi , \p_3 \Phi  \text{ are continuous  on }\R^3\setminus \es ,\\
&\p_1 \Phi (x) \to 0\qquad \text{ as }|x|\to \infty,\\
&\Phi, \p_1 \Phi \in L^2_{\rm loc} (\R^3), \nabla \Phi \in L^2_{\rm loc} (B(0,R)\cap (\R^3\setminus ( \es \cup \widetilde{\es} )))
\end{cases}
\eqne
for every $R>0$. We note that $\nabla \Phi$ is considered as the weak derivative on $\R^3 \setminus (\es \cup \widetilde{\es })$, while $\p_1 \Phi$ is the weak derivative on $\R^3$. Since $\es$ is placed in the uniform shear flow $\mathcal{V}$ we will only consider
\eqnb\label{w_choice}
w\coloneqq -V \sin \alpha .
\eqne
The fourth condition of \eqref{LSP} reflects the fact that the pressure function resulting from placing $\es$ in the flow is localized to $\es$. 
The continuity conditions in \eqref{LSP} are understood in the sense that 
\[
[ \p_j \Phi (x)]_-^+ =0
\]
for each $x\in \widetilde{\es}$, $j\in \{1,3\}$, where we set
\[
[ \p_j \Phi (x)]_-^+  \coloneqq \lim_{\substack{y\to x \\ y_3>0 }} \p_j \Phi (y) - \lim_{\substack{y\to x \\ y_3<0 }} \p_j \Phi (y)
\]
for such $x,j$. Note that $\p_j \Phi (y)$ is well-defined pointwise in $\R^3\setminus (\es \cup \widetilde{\es })$ as a harmonic function.  

The system \eqref{LSP} is known as the \emph{lifting surface problem}, and it was first studied by Multhopp~\cite{multhopp}.

\begin{rem}[The importance of \eqref{LSP}]
The importance of the lifting surface problem is that it is the main  model which can be used to compute lift and aerodynamic drag of a flat plate. Namely, when a solution $\Phi$ to \eqref{LSP} is obtained, we can integrate the difference of  pressure function $P$  (given by \eqref{P_def}) below and above $\es$ to obtain the force normal to the wing,
\eqnb\label{normal_force}
N = -\int_{\es} \left[ P \right]^+_- =  \rho V \cos \alpha  \int_{\es} \left[ \p_1 \Phi  \right]^+_- .
\eqne
The lift force is then the component of $N$ in the direction perpendicular to $\mathcal{V}$,
\eqnb\label{potential_lift}
{L_p} = N \cos \alpha  =   \rho V \cos^2 \alpha  \int_{\es} \left[ \p_1 \Phi  \right]^+_-  ,
\eqne
where the first equality comes from the fact that $-\int_{\es} \left[ P \right]^+_-$ denotes the force acting in the direction perpendicular to $\es$ (so that multiplying by ``$\cos \alpha$'' we obtain the force acting in the direction perpendicular to $\mathcal{V}$).  We will refer to  \eqref{potential_lift} as the \emph{potential lift} of $\es$. 
\end{rem}
\begin{rem}[The leading-edge suction]\label{rem_suction}
We emphasize that \eqref{potential_lift} \emph{is not the total lift} generated by $\es$, even when we  restrict ourselves to the potential flow only (that is relying only on the solution $\Phi$ to \eqref{LSP}). Indeed, the sharp edges of the leading edge of $\es$ result in the additional \emph{leading-edge suction force} 
\eqnb\label{suction_force}
T = \rho V \sin \alpha \int_{\es} \left[ \p_1 \Phi  \right]^+_-  
\eqne
in the plane of the wing, see Fig.~\ref{fig_suction} for a sketch. The appearance of $T$ is caused by the sharp leading edge. To be precise, the stagnation point of the fluid flow occurs away from the leading edge and below the wing (see Fig.~\ref{fig_analogy}(a)), and the fluid particles leaving the stagnation points on the side of the leading edge must, roughly speaking, ``turn around'' the leading edge. By the conservation of total momentum, this will cause the leading edge to be pulled in the direction of the incoming flow. This is the leading edge suction, and the formula \eqref{suction_force} can be proved rigorously in the  case of $2$D flows around a flat airfoil (see \cite[Section~6.5.2]{katz_plotkin}). We only treat $T$ in \eqref{suction_force} above as the 3D analogy. 

\begin{figure}[htbp]
\centering
    \includegraphics[width=9cm]{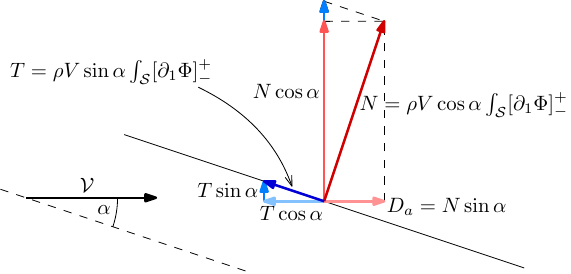}
  \caption{A sketch (in the plane perpendicular to $x_2$) of the role of the suction force $T$ in a sharp leading edge airfoil, in the context of the Kutta-Joukovski theorem. Note that the pressure jump across the airfoil produces the normal force $N$, which contributes to both the lift force $N \cos \alpha$ and the aerodynamic drag force $D_a$. On the other hand, the leading-edge suction $T$ has horizontal component which cancels $D_a$, while its vertical component contributes to $N\cos \alpha$ to give the total lift-force $\rho V \int_{\es} [\p_1 \Phi ]^+_-$. }   \label{fig_suction}
\end{figure}

We now comment how formula \eqref{suction_force} can be understood better in the context of the Kutta-Joukovski theorem and the aerodynamic drag.
Considering the component of $T$ that is parallel to $\mathcal{V}$ we obtain
\[
T \cos \alpha = \rho V \sin \alpha \cos \alpha \int_{\es} \left[ \p_1 \Phi  \right]^+_-  
\]
This is exactly equal to the aerodynamic drag, obtained analogously to \eqref{potential_lift},
\[
D_a = -\sin \alpha\int_{\es} \left[ P \right]^+_- =  \rho V \sin \alpha \cos \alpha  \int_{\es} \left[ \p_1 \Phi  \right]^+_-.
\]
In fact, these two forces cancel each other, as predicted by the d'Alembert computation (see \cite{anderson} for the 2D case).

As for the component of $T$ that is normal to $\mathcal{V}$, including it in the lift force gives the total of
\[
L_p + T\sin \alpha = \rho V  \int_{\es} \left[ \p_1 \Phi  \right]^+_-  
\]
This should be understood as a corollary of the Kutta-Joukovski theorem, by denoting the \emph{effective circulation} by
\eqnb\label{eff_Gamma}
\Gamma \coloneqq \frac{1}b \int_{\es} \left[ \p_1 \Phi  \right]^+_-  ,
\eqne
so that the lift is $L= \rho V \Gamma $ per unit span.  In particular, the jump of $\p_1 \Phi$ across $\es$ should be thought of a vector measure in the direction of $x_2$, representing vorticity. We emphasize that, strictly speaking, the Kutta-Joukovski theorem (see \cite[(3.140)]{anderson} for example) is only a statement for 2D flows, and  the above explanation is only a 3D interpretation.
\end{rem}

\begin{rem}[\eqref{LSP} as an inverse problem]\label{rem_inv}
We note that, as observed by \cite{multhopp}, using the linearized equations \eqref{lin_euler}, the downwash $w(x_1,x_2)$ can be determined from $[ P ]^+_-$  using the double-layer potential theory,
\eqnb\label{DLP}
w(x_1,x_2) = \frac{1}{4\pi \rho V \cos \alpha } \mathrm{f.p.}\int_{\es} \frac{[P(y_1,y_2)]^+_-}{(x_2-y_2)^2 } \left( 1 + \frac{x_1-y_1}{|x-y|} \right) \d y.
\eqne
Indeed, $P$ solves $\Delta P = 0$ in $\R^3 \setminus (\es \cup \widetilde{\es} )$ and is continuous on $\R^3 \setminus \es$. Applying the double-layer potential (see \cite{mclean} for details) to represent $P$ in $\R^3$, deducing from \eqref{P_def} that $\Phi = -\frac{1}{\rho V \cos \alpha } \int_{-\infty }^{x_1} P(y,x_2,x_3) \d y$, and taking $\p_3$ leads to \eqref{DLP}, via a standard computation, see \cite[(15)]{multhopp} and the alternative derivation \cite[Appendix~I]{multhopp} of \eqref{DLP}. 
The ``$\mathrm{f.p.}$'' in \eqref{DLP} denotes the Hadamard finite part. 

Thus, the lifting surface problem \eqref{LSP} can be thought of as an inverse problem to \eqref{DLP}, in which case \eqref{LSP} should be amended by the Kutta-type condition at the trailing edge, such as continuity of $P$ on $\R^3 \setminus \es$, which is stated as a part of the third condition in \eqref{LSP}. 

Moreover, the explicit integral representation \eqref{DLP} of $w$ makes the lifting surface problem \eqref{LSP} susceptible to discretization and numerical solution, which is the main thrust of Multhopp~\cite{multhopp}, see also~\cite{van_spiegel_wouters,lamar}, who refined the method  and the the numerical quadrature. We explain the method in Section~\ref{sec_comp_Kp} in the case when $\es$ is a delta surface.
\end{rem}

\begin{rem}[The jump of the potential in the wake]
We emphasize that $\Phi$ is allowed to be discontinuous across the wake $\widetilde{\es}$, which is an important aspect of the problem.

In fact, is appears to be an issue that is often misunderstood even in the case of 2D airfoils. A good example of this is the d`Alembert paradox, which is a theorem saying that, if $u$ is a potential flow outside of an airfoil (with reasonable decay as $|x|\to \infty$), then the lift generated on the airfoil is $0$ (see \cite[Section~4.13]{acheson} or \cite[Section~8.2]{sverak_notes} for a proof). This paradox is often misunderstood as statements that ``airplanes cannot fly in inviscid fluid'', or ``2D airfoils in inviscid fluid generate zero lift''. These are false statements: the velocity field $u$ in the standard $2$D airfoil theory is potential only locally, and is not potential in the entire outer domain of the airfoil, as the potential admits a jump across the wake (while the velocity itself remains smooth in the entire outer domain), with the jump value equal to the total circulation around the airfoil, giving lift by the Kutta-Joukovski theorem. This is a well-known fact, which was known already to  Helmholtz~\cite{helmholtz} and Lamb~\cite{lamb_1879}, see also chapters 4.9-11 in the book~\cite{acheson} of Acheson for a clear rigorous explanation. The ``paradox'' is really about the drag force, as in practice viscous effects cause skin friction and flow separation which \emph{always} produce finite drag, see  \cite[p.~210]{anderson}. A good explanation was provided by Prandtl in his landmark 1904 paper \cite{prandtl_04}  introducing the  boundary layer theory. 

In the case of 3D airfoil, modelled by a flat surface followed by a linearized wake (as in Fig.~\ref{fig_sketch} below), the problem of the potential jump at the wake is more challenging. For example, in \eqref{LSP} the potential jump can vary in $x_2$, and so, in particular, the velocity field defined by $\nabla \Phi$ in $\R^3 \setminus (\es \cup \widetilde{\es})$ will admit a jump across $\widetilde{\es}$, resulting in the \emph{trailing vortex sheet}. The question of determining the potential jump for any given flat surface $\es$ is a very interesting open problem, as is the lifting surface problem \eqref{LSP} in general. 
\end{rem}

\begin{rem}[known solution methods of \eqref{LSP}]
As mentioned in Remark~\ref{rem_inv}, the lifting surface problem \eqref{LSP} can be approximated using a number of numerical methods. Modern techniques also involve the Vortex Lattice Method (VLM), see \cite[Section~5.5]{anderson} or  Drela \cite{drela} for a modern exposition. 
Moreover, an analytical solution to \eqref{LSP} has only been found in the case when $\es$ is a  disc, in a remarkable technical treatise of Boersma~\cite{boersma}. His analysis in this special case can also be discretized to yield more accurate numerical predictions, see~\cite[Section~5.6]{boersma}. Some consequent analytical work \cite{jordan1, jordan2, kida, haupt_miloh} made some progress in the case when $\es$ is an ellipse. 
\end{rem}

\section{The delta wings and the leading edge suction analogy}\label{sec_polhamus}

In this section we describe Polhamus' \cite{polhamus} analysis of delta wings, which leads to the formula \eqref{polhamus_ansatz}.  

\begin{figure}[htbp]
\centering
    \includegraphics[width=13cm]{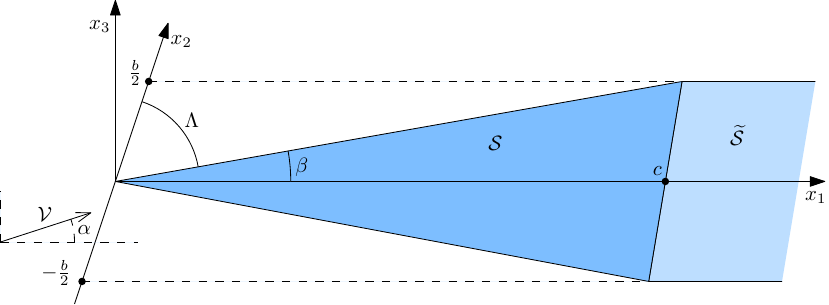}
  \caption{A sketch of the delta wing.}   \label{fig_sketch}
\end{figure}

We will consider a delta wing $\es$ of a fixed area $|\es|>0$, wingspan $b>0$, so that the aspect ratio $A = b^2/|\es| = 4 \tan \beta $ (recall~\eqref{aspect_ratio}), where  $2\beta \in (0, \pi )$ is the angle at the wingtip, and we can write
\eqnb\label{es_def}
 \es\coloneqq \left\lbrace x \colon  |x_2 |  \leq  \frac{b}2 , \frac{x_2}{\tan \beta }  \leq x_1 \leq  \frac{b}{2 \tan \beta } , x_3 =0 \right\rbrace ,
\eqne
see Fig.~\ref{fig_sketch}. We also define the leading edge angle 
\[
\Lambda \coloneqq \frac{\pi}2 - \beta ,
\]
the wing length 
\eqnb\label{c_def}
c \coloneqq \left( \frac{|\es|}{\tan \beta } \right)^{1/2} = \frac{b}{2 \tan \beta } = \frac{b}2 \tan \Lambda = \frac{2 b }A,
\eqne
and we define the wake $\widetilde{\es}$ as in \eqref{wake_def}, see Fig.~\ref{fig_sketch}.

\subsection{The leading edge suction analogy}

Here we describe the arguments of Polhamus~\cite{polhamus} leading to the derivation of \eqref{polhamus_ansatz}.

The potential lift of a delta wing can be computed using the lifting surface problem \eqref{LSP}. Namely, \eqref{normal_force} gives the normal force $N= \rho V \cos \alpha \int_{\es} [\p_1 \Phi ]^+_-$, \eqref{potential_lift} gives the potential lift $L_p = N \cos \alpha$, so that the lift coefficient $C_{L,p}$ due to the potential lift becomes 
\eqnb\label{lift_coef_pot}
C_{L,p} \coloneqq \frac{L_p}{\frac12 \rho V^2 |\es|} = \frac{\cos^2\alpha }{\frac12 V |\es| } \int_{\es} [\p_1 \Phi ]_-^+ 
\eqne
Here the potential the potential $\Phi$ is a solution of \eqref{LSP}--\eqref{w_choice}, and so $\Phi$ is linear with respect to  $w= - V \,\sin \alpha $. Thus there exists a dimensionless constant $K_p$, the potential lift slope, depending only on the aspect ratio $A$, such that
\eqnb\label{Kp_def}
C_{L,p} = K_p \, \sin\alpha \, \cos^2 \alpha ,
\eqne
cf. \cite[(5)]{polhamus}. Polhamus \cite{polhamus} actually states  that  that $K_p$ can be determined from any ``suitable lifting-surface theory'', and that he determines $K_p$ (as a function of $A$; see Fig.~\ref{fig_KpKv} for the result), using the numerical solution of \eqref{LSP} described by Multhopp \cite{multhopp}, with the modifications due to Van Spiegel and Wouters~\cite{van_spiegel_wouters} and Lamar~\cite{lamar}.  As mentioned in Remark~\ref{rem_inv} above, the method originates from discretization of the double-layer potential representation \eqref{DLP}, and we describe this method in detail in Section~\ref{sec_comp_Kp} below. \\

\begin{figure}[htbp]
\centering
    \includegraphics[width=9cm]{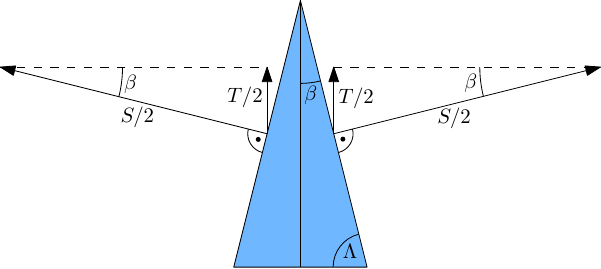}
  \caption{The leading edge suction $T$ as a complement of the suction force $S$. }   \label{fig_suction_force}
\end{figure}

\begin{figure}[htbp]
\centering
    \includegraphics[width=7cm]{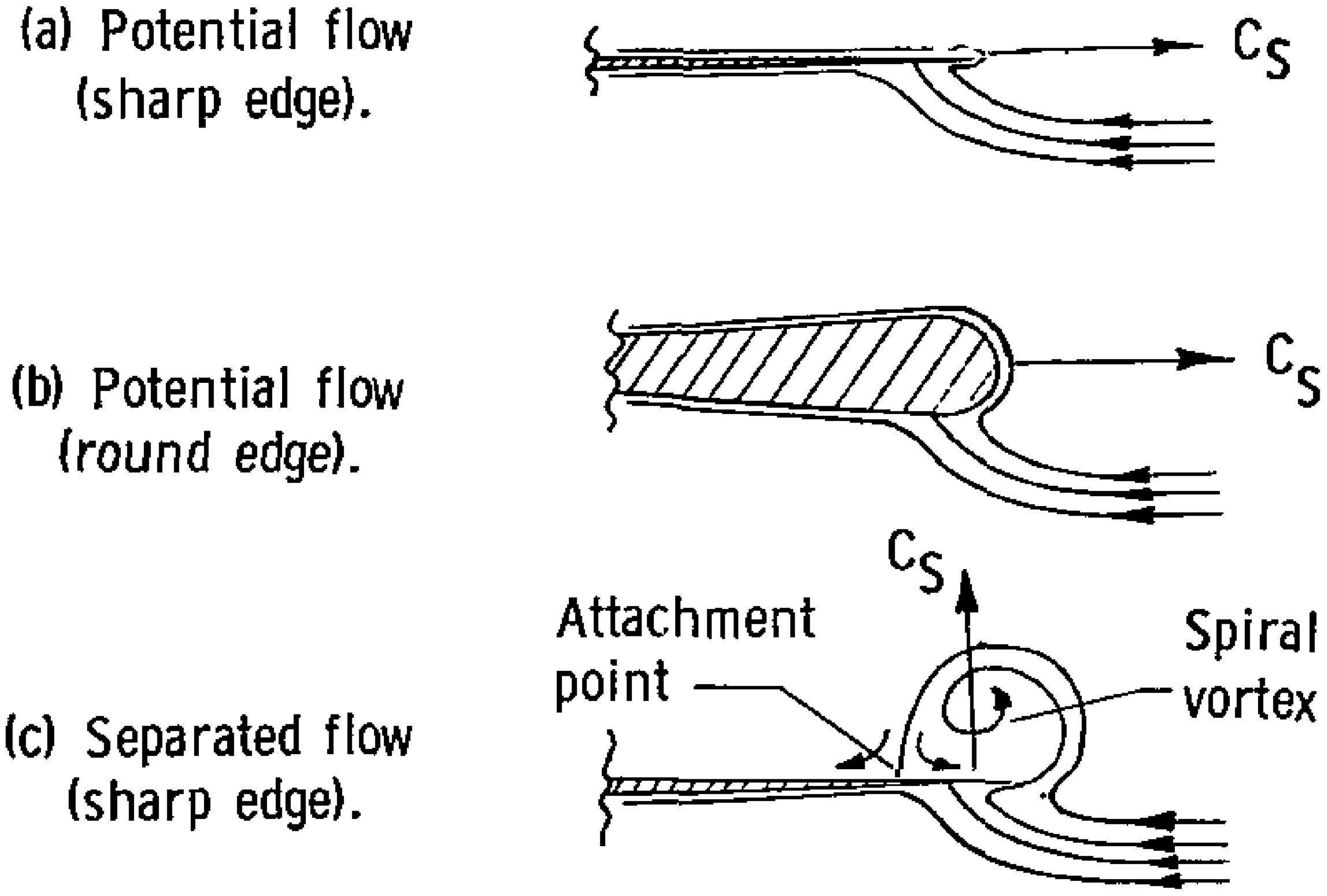}
  \caption{(after \cite[Fig.~7]{polhamus}) An explanation of the leading-edge suction analogy. Here $C_S \coloneqq S/\frac{1}2 \rho V^2 |\es |$ is the suction coefficient.}   \label{fig_analogy}
\end{figure}

The main contribution of Polhamus is concerned with interpretation of the additional force, the \emph{vortex lift} of the total lift force,   originates from the attached vortices, as sketched visible in Fig.~\ref{fig_sketch_vortices}. 
In order to understand this, we first note that we can state the normal force $N$ in terms of an \emph{effective circulation} $\Gamma$ (recall~\eqref{eff_Gamma}), 
\eqnb\label{Gamma_prop}
\Gamma = \frac{1}b  \int_{\es} [\p_1 \Phi ]_-^+  = K_p \frac{|\es|V}{2b } \sin \alpha ,
\eqne
where we used the definition of $K_p$ in the second step, so that the potential lift $L_p$ can be stated in the Kutta-Joukovski form
\[
N  = \rho V \cos \alpha \, \Gamma b ,
\]
cf. \cite[(2)]{polhamus}. The effective circulation $\Gamma $ makes it easy to compute the leading edge suction of the airfoil. Indeed, it suffices to replace the component $V \cos \alpha$ that is tangential to $\es$ by its normal component $V \sin \alpha$. This gives the leading edge-suction force in the $-x_1$ direction
\[
\rho V \sin \alpha \, \Gamma b
\]
as pointed out in Remark  and Fig.~\ref{fig_suction}, this suction force also contributes to the total lift. However, Polhamus's approach assumes a different treatment. First, in order to obtain the total suction force $T$, he takes into account the downwash velocity $w_i$ (in addition to the normal component $V\sin \alpha$ of the free stream), which arises from the trailing vortex sheet (see below for a discussion of $w_i$), which gives that
\eqnb\label{T_def_analogy}
T = \rho \Gamma b  \left( V \sin \alpha - w_i \right).
\eqne
Having found the leading-edge suction force $T$, the main idea of Polhamus is a modelling assumption consisting of two aspects:
\begin{enumerate}
\item Suppose that $T$ is a component of a suction force $S$ perpendicular to the leading edges and in the plane of the wing, see Fig.~\ref{fig_suction_force}. Using this assumption we write
\[
T =  S \sin \beta = S \cos \Lambda.
\]
\item Suppose that, $S$, rotated into the direction perpendicular to the wing plane, can be treated as the vortex lift. This modelling assumption is based on the similarity between flow streamlines of a 2D past a 2D airfoil and the observed flow and its reattachment in the case of delta wings of small aspect ratio. 
\end{enumerate}
The above two ideas are the \emph{leading edge suction analogy}. Point (2) can be better understood through the sketch in Fig.~\ref{fig_analogy} (see also Fig.~\ref{fig_sketch_vortices}), which shows that fluid particles first move away from the lateral edges of the wing in the vertical direction (that is $x_3$) and then reattach to the wing on the other side of the attached vortices. Thus, analogously to the explanation of the leading-edge suction (in Remark~\ref{rem_suction}), the conservation to total momentum implies that a suction force in the vertical direction will be exerted on the wing. Polhamus~\cite[p.~8]{polhamus} explains  that the ``total force is associated with the pressures required to maintain the equilibrium of the flow over the separated spiral vortex is essentially the same as the leading-edge suction force associated with the
leading-edge pressures required to maintain attached flow around a large leading-edge
radius'' and refers to Fig.~\ref{fig_analogy}.

The total normal force is thus 
\eqnb\label{total_lift_polhamus}
N + \frac{T}{\cos \Lambda } = \rho b \Gamma V \cos \alpha  + \rho \Gamma b \frac{V \sin\alpha - w_i }{\cos \Lambda },
\eqne
where we used \eqref{lift_coef_pot}, \eqref{Gamma_prop} and \eqref{T_def_analogy}.
As for the downwash velocity $w_i$, Polhamus points out that it must be proportional to $\Gamma$, which in turn is proportional to $V\sin \alpha$  (by~\eqref{Gamma_prop}), and states that 
\[
w_i = K_p K_i V \sin \alpha,
\]
where $K_i \coloneqq \frac{\p C_{D_i}}{\p C_L^2}$, cf. \cite[(9)]{polhamus}. Substituting this into \eqref{total_lift_polhamus} gives 
\[
N + \frac{T}{\cos \Lambda } = \rho b \Gamma V \cos \alpha  + \rho \Gamma b \frac{1 - K_p K_i }{\cos \Lambda } V \sin\alpha .
\]
Substituting $\Gamma$ from \eqref{Gamma_prop} and taking the vertical component we thus obtain the total lift force
\[
L = N\cos \alpha + \frac{T\cos \alpha }{\cos \Lambda } = \frac{1}2 \rho |\es | K_p   V^2 \sin \alpha  \cos^2 \alpha  +  \frac12 \rho |\es |  \frac{K_p - K_p^2 K_i }{\cos \Lambda } V^2 \sin^2 \alpha \cos \alpha  =: L_p + L_v
\]
That is the total coefficient
\eqnb\label{total_lift_polhamus1}
C_L \coloneqq \frac{L}{\frac12 \rho V^2 |\es |} = K_p \sin \alpha  \cos^2 \alpha  + K_v \sin^2 \alpha  \cos \alpha 
\eqne
as stated in \eqref{polhamus_ansatz}.

\subsection{Explanation of the Polhamus approach}\label{sec_polhamus_explanation}

Here we discuss a more precise interpretation of the above derivation regarding the computation of constants $K_p$ and $K_i$.

As for $K_p$, its value may be computed by attempting a numerical solution of the lifting surface problem \eqref{LSP}, which is what he meant by a ``suitable lifting-surface theory''. We outline this in the following Section~\ref{sec_Kp}.

The computation of $w_i$ and $K_i$, which is not described in the original work \cite{polhamus} in any detail, is more nuanced. 
The constant $K_i = \frac{\p C_{D_i}}{\p C_L^2}$ is a standard notation for computing  the ratio of the induced drag $D_i$ to lift, and it is considered at angle of attack $\alpha =0$. In particular, the assumption that $w_i$ is proportional to $V \sin \alpha$, implicitly assumes that $w_i$ is constant, i.e. that an ``effective downwash velocity'' is meant here. Actually, for delta wing of small aspect ratio $A$, the downwash velocity \emph{is} constant, which only follows as consequence of Proposition~\ref{P1} and the computation~\eqref{optimal_comp} for the 2D model. 
The coefficient $K_i$ can only be computed analytically using the 
Munk lifting-line model of the trailing vortex sheet (which, for completeness, we introduce in detail in Appendix~\ref{app_munk}). It is thus reasonable to assume that, after computing $K_p$, Polhamus 
\eqnb\label{pol_change}
\text{\emph{changes the modelling assumption of the delta wing into the lifting-line model}}
\eqne
and uses the computation of $K_p$ to determine (numerically) the circulation distribution $\Gamma (x_2)$, that is the circulation (in the $x_1-x_3$ plane) at spanwise location $x_2 \in (-b/2,b/2)$. 

Knowing $\Gamma (x_2)$, one can compute the downwash velocity $w_i = w [ \Gamma ](x_2)$ at each spanwise location $x_2$ using the Munk lifting-line model (see~\eqref{downwash_formula} for details). The computation shows that 
\eqnb\label{polh_wi_const}
w_i (x_2)\quad  \text{\emph{appears to be constant},}
\eqne
which is perhaps the reason why Polhamus~\cite[(7)--(9)]{polhamus} does not discuss its spatial dependence.

After computing the downwash $w[\Gamma ]$, one easily computes the induced drag
\[
D_i = \rho \int_{-\frac{b}2}^{\frac{b}2} \Gamma (x_2 ) w [\Gamma ] (x_2 ) \d x_2
\]
(see~\eqref{ind_drag} for details). Moreover, note that, by~\eqref{total_lift_polhamus1}, the vortex lift is of higher order in $\alpha$, so that in order to compute the  ratio $K_i$, it suffices to consider the lift coefficient $C_{L,p}$ of the potential lift only. In that case the potential lift $L_p$ can be computed simply by integrating $\Gamma (x_2)$, $L_p = \rho V \int_{-\frac{b}2}^{\frac{b}2} \Gamma $ (which is equivalent to \eqref{lift_coef_pot}). In that case we obtain
\[
K_i = \frac{1}{\pi A e},
\]
where $e \coloneqq a_0^2 / \left( \sum_{k\geq 0 } (2k+1)a_{2k}^2 \right)$ denotes the efficiency parameter, and the $a_{2k}$'s are the coefficients in the expansion \eqref{Gamma_exp} of $\Gamma $ in Chebyshev polynomials, see Appendix~\ref{app_munk} for details.  Clearly $e\leq 1$ and $e =1$ if and only if $a_{2k}=0$ for all $k\geq 1$. 

It turns out that the delta wing for small aspect ratio becomes maximally efficient as $A\to 0$, which is a consequence of Proposition~\ref{P1}, see the computation~\eqref{optimal_comp}, which shows that a $2$D approximation considered here gives the elliptic circulation distribution, that is $a_{2k}=0$ for all $k\geq 1$. This shows that we can take $K_i = 1/\pi A$, so that 
\eqnb\label{Kv_delta}
K_v = \left( K_p - \frac{K_p^2 }{\pi A} \right) \frac{1}{\cos \Lambda },
\eqne
which is presummably what Polhamus meant when obtaining the result in Fig.~\ref{fig_KpKv}.

Further comments about the leading-edge suction analogy: The potential lift is not the total lift obtained from the lifting surface theory (by a different treatment of the leading-edge suction). In particular, the aerodynamic drag is neglected in this approach.

It is remarkable that the lifting line works for small aspect ratio.

\subsection{Detailed explanation}\label{sec_Kp}

Here we discuss the algorithm for computing $K_p$ in \eqref{Kp_def}, which was implicitly used by Polhamus~\cite{polhamus}, and how to determine the circulation distribution  $\Gamma (x_2)$ (as mentioned below \eqref{pol_change}).

The main idea of Polhamus \cite{polhamus} is to use the observation \eqref{DLP} to discretize the problem and find discretized pressure $p$ from $w$ (note $w$ is constant so it is easy to discretize). Then use the discretized pressure to find discretized circulation for each spanwise location.

\subsubsection{Computing $K_p$}\label{sec_comp_Kp}

\begin{figure}[htbp]
\centering
    \includegraphics[width=13cm]{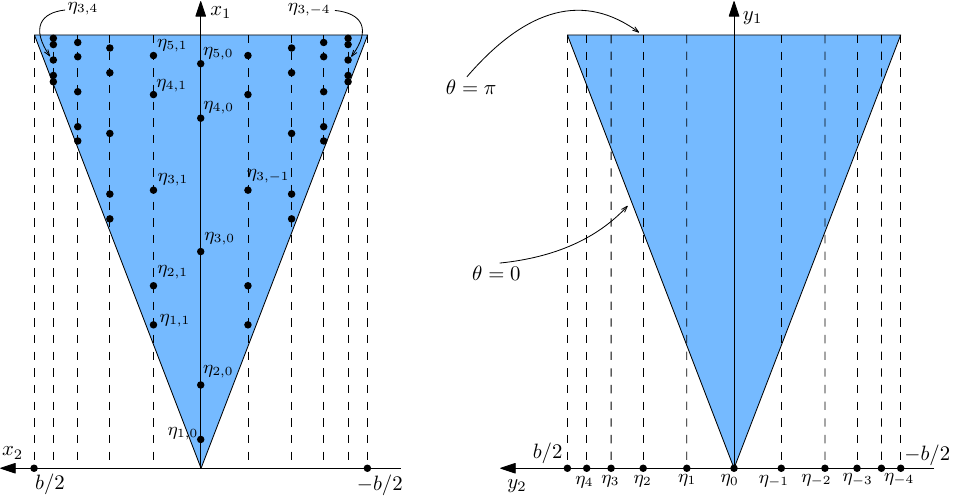}
  \caption{A sketch of the choice of the $\eta_{s,l}$'s and $\eta_m$'s. Here $n=4$, $N=5$. For clarity, the labels of some of the $\eta_{s,l}$'s were omitted in the sketch.}   \label{fig_points}
\end{figure}

For each $y_2 \in (-b/2,b/2)$ we parametrize 
\[
y_1 \in \left( \frac{y_2 }{\tan\beta} ,  \frac{b/2}{\tan\beta} \right) 
\]
using $\theta \in (0,\pi)$, 
\eqnb\label{y1_param}
y_1  \coloneqq  \left( \frac{b/2+y_2 }{2 {\tan\beta}} - \frac{b/2-y_2 }{2 {\tan\beta} } \cos \theta \right) ,
\eqne
so that $\theta =0$ corresponds to the leading edge, and $\theta = \pi $ corresponds to the trailing edge, see Fig.~\ref{fig_points}.
Using this parametrization we can use a finite Fourier Sine expansion in $y_1$ to approximate $p$ in $y_1$.

As for the $y_2$ variable, 
One can consider piecewise polynomial approximation of the $y_2$ dependence, as is done in \cite{lamar}, but here, for the sake of simplicity, we only consider piecewise affine approximation. To this end, given an odd integer $2n+1$ and $m \in \left\lbrace -n, \ldots , 0 ,\ldots , n \right\rbrace$, we set 
\[
\eta_m \coloneqq \frac{b}2 \sin \left( \frac{m \pi }{2n+2} \right) 
\]
(see Fig.~\ref{fig_points}) and we let $\chi_m \in C_c (\R)$ be piecewise affine and such that 
\[
\chi_k (\eta_m) = \delta_{mk}
\]
for all $m,k \in \{ -n ,\ldots , n \}$. 
We will abuse the notation slightly by writing $p(\theta , \eta ) \equiv p (y_1,y_2)$ for $(y_1,y_2,0) \in \es$. 

We will approximate $p$ by
\eqnb\label{p_change_vars}
p(\theta, y_2) =  \sum_{m= -n}^n \sum_{j=0 }^{N-1} \frac{{b/2}- y_2}{\tan\beta} p_{m,j} \chi_m (y_2) \phi_j (\theta ), 
\eqne
where 
\[
\phi_j (\theta ) = \begin{cases}
\cot \frac{\theta }2 \hspace{1cm} & j=0,\\
\sin (j \theta ) &j\geq 1
\end{cases}
\]
denote the chordwise Fourier Sine modes for $j\geq 0$. The lowest order mode  $\cot \frac{\theta }2 $  captures the singularity (which is expected to be of $x^{-1}$ order due to the sharp edge of the wing). As for the downwash $w(x_1,x_2)$, we will apply a similar approximation, except that we will use piecewise affine chordwise approximation. Namely, for each $l \in \left\lbrace -n, \ldots , 0 ,\ldots , n \right\rbrace$, $s \in \{ 1, \ldots , N \}$ we set
\[
\eta_{s,l}  \coloneqq \frac{\eta_l}{\tan\beta} + \frac{{b/2}-\eta_l}{2{\tan\beta} } \left(1- \cos \frac{2\pi s }{2N+1} \right)
\]
we let $\chi_{s,l} \in C_c (\R^2)$ be a piecewise affine function such that $\chi_{s,l}(\eta_{s',l'})= \delta_{ss'}\delta_{ll'}$,
and we approximate $w(x_1,x_2)$ by 
\eqnb\label{w_approx}
w (x_1,x_2) = \sum_{l=-n}^n \sum_{s=1}^{N} w_{sl} \chi_{s,l} (x_1,x_2) 
\eqne
Under such approximation the double-layer potential formula \eqref{DLP},
\[
4\pi \rho V \cos \alpha \, w (x_1,x_2) = \mathrm{f.p.} \int_{\es} \frac{[P(y_1,y_2)]^+_-}{(x_2-y_2)^2 } \left( 1 + \frac{x_1-y_1}{|x-y|} \right) \d y,
\]
becomes
\eqnb\label{lin_sys}
\begin{split}
4\pi \rho &V \cos \alpha \,\sum_{l=-n}^n \sum_{s=1}^{N} w_{sl} \chi_{s,l} (x_1,x_2) \\
&\hspace{1cm} \sim   \sum_{m= -n}^n \sum_{j=0 }^{N-1} p_{m,j}   \mathrm{f.p.} \int_{-{b/2}}^{b/2} \int_{0}^\pi \frac{\chi_m (y_2) \phi_j (\theta ) }{(x_2-y_2)^2 } \left( 1 + \frac{x_1-y_1(\theta )}{|x-(y_1(\theta ), y_2 )|} \right)\sin \theta \, \d \theta\, \d y_2,
\end{split}
\eqne
where we also recalled \eqref{y1_param} to the change  variable $y_1\mapsto \theta$. For each $l$, $s$ we now take $(x_1,x_2) = \eta_{s,l}$ to get
\eqnb\label{lin_sys1}
4\pi \rho V \cos \alpha \,\, w_{s,l} =  \sum_{m= -n}^n \sum_{j=0 }^{N-1} L_{(s,l),(m,j)} p_{m,j} , 
\eqne
where $L_{(s,l),(m,j)}$ denotes the double integral on the right-hand side of \eqref{lin_sys} evaluated at $(x_1,x_2)=\eta_{s,l}$. (In \cite{lamar}, a more sophisticated quadrature is used to find  $L_{(s,l),(m,j)}$.) In vector notation $w \equiv \{ w_{s,l} \} \in \R^{N(2n+1)}$, $p \equiv \{ p_{m,j} \} \in \R^{N(2n+1)}$, $L = \{L_{(s,l),(m,j)} \} \in \R^{N(2n+1) \times N(2n+1)}$ we thus obtain a linear system
\eqnb\label{lin_sys_vec}
4\pi \rho V \cos \alpha \,\, w =  L\,p.
\eqne
Inverting this linear system we can thus find $p =\frac{1}{ 4\pi \rho V \cos \alpha } L^{-1} w$, where $w_{s,l} = - V \sin \alpha$ for all $s,l$, due to \eqref{w_choice}.

Having found $p_{m,j}$ for all $m,j$, we can compute lift by approximating \eqref{potential_lift},
\eqnb\label{lift_temp1}
\begin{split}
L_p &= -\cos \alpha \int_{\es} [p]^+_- = - \cos \alpha \int_{-{b/2}}^{b/2} \int_0^\pi [p( \theta , y_2 )]^+_- \sin \theta \, \d\theta \,\d y_2 \\
&\sim - \cos \alpha \sum_{m= -n}^n \sum_{j=0 }^{N-1} p_{m,j} \int_{-{b/2}}^{b/2}  \int_0^\pi \chi_m (y_2 ) \phi_j (\theta ) \sin \theta \, \d\theta \,\d y_2 \\
&= - \cos \alpha \sum_{m=-n}^n   \left( p_{m,0}\underbrace{\int_0^\pi \cot\frac{\theta }2 \sin\theta \d \theta }_{=\pi }+ p_{m,1}  \underbrace{ \int_0^\pi  \sin^2 \theta \d \theta}_{=\pi/2} \right)  \int_{-{b/2}}^{b/2} \chi_m\\
&=\frac12 \rho V^2 |\es| \cdot   \cos^2 \alpha \sin \alpha    \sum_{m=-n}^n  \frac{8\pi^2}{|\es|} \left( (L^{-1} h )_{m,0}  + \frac{(L^{-1} h )_{m,1}}2  \right)  \int_{-{b/2}}^{b/2} \chi_m
\end{split} 
\eqne
where $h\coloneqq (1,1,\ldots , 1 ) \in \R^{N(2n+1)}$ and, in the third line, we noted that $\int_0^\pi \sin (j\theta ) \sin \theta \d \theta =0$ for $j\geq 2$. 
Since $L_p = \frac12 \rho V^2 {|\es|} \cdot K_p \cos^2\alpha\, \sin \alpha $ (recall~\eqref{lift_coef_pot}--\eqref{Kp_def}), we thus obtain that
\[
K_p \sim \sum_{m=-n}^n  \frac{8\pi^2}{|\es|} \left( (L^{-1} h )_{m,0}  + \frac{(L^{-1} h )_{m,1}}2  \right)  \int_{-{b/2}}^{b/2} \chi_m.
\]
Note that the definition of $L$ (below~\eqref{lin_sys1}) shows that $L$ has units $[1/\mathrm{length} ]$, so that the right-hand side above is dimensionless. 

\subsubsection{Finding  $\Gamma (x_2)$} Having found $p_{m,j}$ in \eqref{lin_sys_vec} we observe that, for each $m$, $p_{m,0},p_{m,1}$ determine the lift contribution at the spanwise location $\eta_m$. Namely, we see from \eqref{lift_temp1} that the contribution to the normal force $N$ per unit span at $y_2 \sim \eta_m$ is 
\[
- \left( \pi p_{m,0} + \frac{\pi}2 p_{m,1} \right) 
\]
Interpreting this lift contribution as the lift generated by a circulation $\Gamma (\eta_m)$ in the uniform flow $V\cos \alpha$ (as in \eqref{normal_force}--\eqref{potential_lift})  via the Kutta-Joukovski theorem, we get that
\[
\rho V \cos\alpha \Gamma (\eta_m ) = \pi p_{m,0} + \frac{\pi}2 p_{m,1},
\]
Consequently, we can approximate 
\[
\Gamma (y_2 ) = \frac{1}{\rho V \cos \alpha } \sum_{m=-n }^n \left( \pi p_{m,0} + \frac{\pi}2 p_{m,1} \right)  \chi_m (y_2) .
\]
As mentioned above, thanks to the  circulation distribution $\Gamma (y_2)$ one can follow the modelling assumption \eqref{pol_change} to replace our model \eqref{p_change_vars}--\eqref{w_approx} of the delta wing  by the lifting line passing through the center of mass, with circulation distribution $\Gamma (y_2)$.

\section{A 2D approach to the problem of delta wing potential flow}\label{sec_2d_approach}

Here we describe a 2D method of approximating the solution $\Phi$ to the lifting surface problem \eqref{LSP} in the case of the delta wing~\eqref{es_def} of small aspect ratio $A$. 

We set  
\eqnb\label{gamma_def}
\gamma \coloneqq  \tan \beta = \cot \Lambda =\frac{A}4
\eqne
for brevity (recall the sketch in Fig.~\ref{fig_sketch}). For $x_1 \in (0,c)$ we let $\phi (x_1,x_2,x_3)$ be the solution of the 2D problem (in $(x_2,x_3)$) 
\eqnb\label{2d_pot}
\begin{split}
\Delta_{x_2,x_3} \phi &=0 \hspace{2cm} \text{ in } \R^2 \setminus \{ (x_2,x_3) \colon |x_2|\leq \gamma x_1 , x_3=0 \},\\
\p_3 \phi &= -V \sin \alpha\hspace{1cm} \text{ for  } x_3=0, |x_2|\leq \gamma x_1 ,\\
|\phi (x) | &\to 0 \hspace{1.9cm}\text{ as } |x_3|\to \infty 
\end{split}
\eqne
given by 
\eqnb\label{phi_def}
\phi (x) \coloneqq -V \sin \alpha \left(x_3-  \im \left( \sqrt{ x_2+ix_3-\gamma x_1 }\sqrt{ x_2+ix_3+\gamma x_1  } \right) \right).
\eqne
We also extend $\phi$ by $0$ for $x_1 <0$ and by $\phi (c^-,x_2,x_3)$ for $x_1 >c$.
That \eqref{phi_def} satisfies \eqref{2d_pot} is clear by noting that  $-i V\sin \alpha \sqrt{ x_2+ix_3-\gamma x_1 }\sqrt{ x_2+ix_3+\gamma x_1  } $ is a complex potential (in the complex plane of $x_2+ix_3$) of the potential flow (in the $x_3$ direction). In other words, \eqref{phi_def}  is a solution to a 2D version of the lifting surface problem \eqref{LSP} in the $(x_2,x_3)$ plane.  The explicit formula \eqref{phi_def} is sometimes referred to as the slender wing theory by Jones \cite{jones}\footnote{Although the exact formula \eqref{phi_def} does not appear explicitly in \cite{jones}.}, and it can be derived from the complex potential of a uniform shear flow, via two conformal transformations, see Appendix~\ref{app_2dpot} for details, where we also explain why we write the product of square roots in \eqref{phi_def}, instead of $\sqrt{(x_2+ix_3)^2 -\gamma^2 x_1^2 }$. \\

As for the regularity of $\phi$, a direct computation shows that 
\eqnb\label{regularity_of_phi}
\phi,\p_1 \phi  \in L^2_{\rm loc } (\R^3 ), \text{ and } \phi \in H^{3/2-\varepsilon } (B(0,R)\cap (\R^3 \setminus (\es \cup \widetilde{\es }) ) ) \text{ for all } R,\varepsilon >0.
\eqne
 Moreover, $\phi \not \in L^2 (\R^3)$ and $\nabla \phi \not \in L^2 (B(0,R)\cap (\R^3 \setminus (\es \cup \widetilde{\es }) ) )$ due to a lack of the decay and slow decay (respectively) as $|x|\to \infty$. Furthermore,  $\phi \not \in H^1 (\R^3 \setminus \es )$ due to the jump of $\phi$ across $\widetilde{\es}$ (so that $\p_3 \phi$ is not even an $L^1_{\rm loc } (\R^3 \setminus \es )$  function). 

Nevertheless, we expect that for $\gamma \ll 1$  \eqref{phi_def} provides a good approximation of the solution $\Phi$ to the lifting surface problem around $\es$. In other words, we see that $\Delta \phi \ne 0$, but we nevertheless take 
\[
u \coloneqq \nabla \phi.
\]
We first note that, on the airfoil plane (i.e. for $x_3=0$) we have 
\eqnb\label{phi_on_line}
\phi (x_1,x_2, \pm 0) =  \pm V \sin \alpha  \sqrt{\gamma^2 x_1^2 - x_2^2 },
\eqne
where the plus and minus signs correspond to the upper and lower surfaces, respectively. Indeed, using the convention that $z\mapsto \sqrt{z}$ admits the branch cut along the negative real axis, we see that, among the two square roots in \eqref{phi_def}, only $\sqrt{x_2+ix_3 -\gamma x_1}$ crosses the branch cut at $x_3 =0$, so that
\[
\lim_{x_3 \to 0^{\pm}} \sqrt{x_2+ix_3 -\gamma x_1} = \pm i | x_2 - \gamma x_1|^{1/2},
\]
which implies \eqref{phi_on_line}. Taking $\p_1$ of \eqref{phi_on_line} we see that 
\eqnb\label{u1_jump}
u_1 (x_1,x_2,\pm 0) =  \pm {V\sin \alpha } \frac{\gamma^2 x_1 }{\sqrt{ \gamma^2 x_1^2- x_2^2 }},
\eqne
which is a formula obtained in \cite[(5.64)]{houghton}. We can now use Bernoulli's law to compute the lift. Namely, we suppose that $\mathcal{V}+u$ satisfies the 3D Euler equations  and we neglect the quadratic terms in $u$ as follows
\eqnb\label{temp0}
\begin{split}
-\frac{1}{\rho} \nabla p &= (\uu+u ) \cdot \nabla (\uu +u ) = \frac12 \nabla |\uu+u |^2 - (\uu+u )\times \underbrace{ \mathrm{curl} (\uu + u )}_{=0}\\
&= \frac12 \nabla |(V\cos \alpha + u_1 , u_2 , V\sin \alpha + u_3 ) |^2 \\
&= \nabla \left( \frac{V^2 \cos^2 \alpha }2  + u_1 V\cos \alpha + \frac{u_1^2}2 + u_2^2 + \frac{(V \sin \alpha + u_3)^2}2 \right) 
\end{split}
\eqne
for $x\in \left( (0,c) \times \R^2 \right) \setminus \es$
where we recalled \eqref{uu_def} in the third equality. This gives the Bernoulii law
\[
p = - \rho \left( \frac{V^2 \cos^2 \alpha }2  + u_1 V\cos \alpha + \frac{u_1^2}2 + u_2^2 + \frac{(V \sin \alpha + u_3)^2}2 \right)  + C
\]
for such $x$, where $C \in \R$ is an arbitrary constant.  We now neglect the quadratic terms $u_1^2/2$, $u_2^2$ and observe that $V\sin \alpha + u_3 =0$ for $x_3=0$ (by~\eqref{2d_pot}) to obtain that 
\eqnb\label{pressure_jump_2D}
\begin{split}
-[ p ]^+_- &\coloneqq  p(x_1,x_2,0^-) -p(x_1,x_2,0^-)  = \rho V \cos \alpha  [ u_1 ]^+_-  = 2 \rho V^2 \sin \alpha \cos \alpha \frac{\gamma^2 x_1 }{\sqrt{\gamma^2 x_1^2  - x_2^2 }}
\end{split}
\eqne
for $(x_1,x_2)\in \es$. 
Integrating $-[ p ]^+_-$ over $\es$ we thus obtain the force acting on the wing in the direction perpendicular to the wing,
\[
N^{(2D)} = - \int_{\es} [p]_-^+.
\] 
This is  analogous to~\eqref{normal_force}.As in~\eqref{potential_lift}, we obtain the potential  lift by  multiplying by $\cos \alpha$,  
\[
\begin{split}
L_p^{(2D)} &= - \cos \alpha \int_{\es} [p]_-^+ = \cos \alpha \int_0^c \int_{-\gamma x_1 }^{\gamma x_1}   [p]^+_- \d x_2\, \d x_1 \\
&=  \rho V^2 \gamma^2 \sin \alpha \cos^2 \alpha \int_0^c \int_{-\gamma x_1}^{\gamma x_1}    \frac{x_1\, \d x_2\, \d x_1}{\sqrt{\gamma^2x_1^2 - x_2^2 }}   \\
&=\frac{\rho}2 V^2 \gamma^2  \sin \alpha\cos^2 \alpha \int_0^{c} x_1 \int_{-1}^1 \frac{\d \xi }{\sqrt{1-\xi^2}} \d x_1  \\
&=2 \rho \pi \tan^2 \beta \, V^2 \sin \alpha\cos^2 \alpha \int_0^{c} x_1 \d x_1\\
& = c^2 \rho \pi \tan^2 \beta \, V^2 \sin \alpha\cos^2 \alpha = \pi  \rho v^2 {|\es|}  \gamma  \sin \alpha \cos^2 \alpha \\
&= \frac12 \rho V^2 {|\es|} \cdot \left( \frac{\pi}2 A \right)   \sin \alpha \cos^2 \alpha 
\end{split}
\]
where we recalled \eqref{c_def} that ${|\es|}=  c^2 \gamma$. Note that, in comparison with \eqref{polhamus_ansatz}, this means that
\eqnb\label{Kp_2d}
K_p^{(2D)} = \frac\pi{2} A,
\eqne
that is the line tangent to Fig.~\ref{fig_KpKv}(left) at the origin. \\

Having found the pressure jump distribution \eqref{pressure_jump_2D} on $\es$, we can compute $\Gamma (x_2)$ using Polhamus modelling assumption. By \eqref{pressure_jump_2D}, 
\[
[ \p_1 \phi ]_-^+ = [ u_1 ]_-^+ = 2 V \sin \alpha \frac{\gamma^2 x_1 }{\sqrt{\gamma^2 x_1^2 -x_2^2 }},
\]
so that the circulation at a fixed spanwise location $x_2$ becomes
\eqnb\label{optimal_comp}
\begin{split}
\Gamma (x_2) &= \int_{\frac{x_2}\gamma }^{\frac{b}{2\gamma }} [ u_1 ]_-^+ \,\d x_1 = 2 V \sin \alpha \int_{\frac{x_2}\gamma }^{\frac{b}{2\gamma }} \frac{\gamma^2 x_1 }{\sqrt{\gamma^2 x_1^2 -x_2^2 }} \,\d x_1 \\
&= 2 V x_2 \sin \alpha \int_{1}^{\frac{b}{2x_2 }} \frac{y }{\sqrt{y^2 -1 }} \d y = \frac{ V  \sin \alpha }b \sqrt{1- \left( \frac2b x_2 \right) ^2 } .
\end{split}
\eqne

This shows that, treating the 2D approximation \eqref{phi_def} of the solution to the lifting-surface problem \eqref{LSP} and changing the modelling assumption to Munk's lifting line model (see Appendix~\ref{app_munk}), the spanwise circulation distribution minimizes the drift-to-lift ratio $K_i$ (see \eqref{Ki_comp}--\eqref{Gamma_optimal}), verifying the comments above \eqref{Kv_delta}. In particular, the induced velocity 
\[
w_i=w[\Gamma ](x_2) = \frac{V\sin\alpha }{2b^2},
\] 
see~\eqref{wGamma_optimal} for details. We can thus treat the induced velocity $w_i$ as constant throughout $\es$,  verifying the observation~\eqref{polh_wi_const} above. 

\begin{rem}[Comments on the 2D potential flow derivation]
Note that the lift does not come from the 2D potential flow $\phi$ itself, but rather from its growth in $x_1$, via $u_1$.

Two main unresolved issues of this $2$D approach. One is that, as pointed out above, even though we used in \eqref{temp0} that $u$ solves the linearized 3D Euler equations, it does not. 

Furthermore the issue of the trailing edge is neglected. Nevertheless the result gives the correct prediction as $\beta \ll 1$.
\end{rem}

\subsection{Justification of 2D approach for small $\Lambda$}\label{sec_pf_P1}

Here we estimate $\phi - \Phi$, that is we prove Proposition~\ref{P1}. For this, let us first consider the problem
\eqnb\label{elliptic_p}
\begin{cases}
-\Delta g = f  &\text{on } \R^3 \setminus (\es \cup \widetilde{\es}) \\
 \p_3 g = 0 &\text{on } \es \\
\p_1 g, \p_3 g \, \text{ are continuous on }\R^3 \setminus \es &\\
g, \p_1 g \in L^2_{\rm loc } (\R^3 )
\end{cases}
\eqne
and we suppose that 
\eqnb\label{the_elliptic_issue}
\text{ Any solution of \eqref{elliptic_p} is bounded }  H^{-1/2}(\R^3 \setminus (\es \cup \widetilde{\es})) \ni f\mapsto g \in H^{3/2} (\R^3 \setminus (\es \cup \widetilde{\es}).
\eqne 

\begin{rem}[about the elliptic regularity~\eqref{the_elliptic_issue}]\label{rem_elliptic}
We emphasize that \eqref{the_elliptic_issue} is not clear. If \eqref{elliptic_p} included extra assumptions that $\p_3 g =0$ on $\widetilde{\es }$ and $\nabla g \in L^2 (\R^3 \setminus (\es \cup \widetilde{\es} ))$ then elliptic regularity $H^{-1 } \ni f\mapsto g \in H^{1} $ could be proved by multiplying the PDE $-\Delta g = f $ by $g\chi$, where $\chi \in C_c^\infty (\R^3) $ is a cutoff function, integrating by parts and applying the Dominated and Monotone Convergence Theorems. The main point of the extra assumption of $\p_3 g =0$ on $\widetilde{\es}$ is to guarantee that the jump of $g$ across $\widetilde{\es}$ does not get in the way of such a proof. Note that the jump $[g]^+_-$ is a function of $x_2$ only, and that $\int_{-b/2}^{b/2} [ g ]^+_- \d x_2= \int_{\es } [\p_1 g ]^+_-$, which could be expected to vanish as the normal force \eqref{normal_force}. Either way, any elliptic theory of \eqref{elliptic_p} relies on a better understanding of the lifting-surface problem~\eqref{LSP}, which is open. It would be particularly interesting to understand how  the discontinuity of the potential $\Phi$ across the wake $\widetilde{\es}$ is related to the shape of $\es$.
\end{rem}

As discussed in the introduction, we pretend that \eqref{the_elliptic_issue} holds, and note that $\phi - \Phi $ satisfies
\[
\Delta ( \phi - \Phi ) = \p_{11} \phi
\]
in $\R^3 \setminus (\es \cup \widetilde{\es} )$, as well the remaining conditions of \eqref{elliptic_p}. Note that  
\[
\p_1 \phi =  - V \sin \alpha \, \gamma^2 x_1 \, \im  \left( (x_2+ix_3)^2 - \gamma^2 x_1 \right)^{-\frac12} = - \gamma^2 x_1 W (\gamma x_1 , x_2,x_3 )
\]
for $x_1 \in (0,c )$, and $\p_1 \phi =0$ for $ x_1 \in (-\infty , 0) \cup (c,\infty )$,
where 
\[
W ( r , x_2,x_3) \coloneqq  \im \frac{V \sin \alpha }{\sqrt{(x_2+ix_3)^2 - r^2 }}.
\]
Moreover, for $x_1 \in (0,c)$,
\[
\p_{11} \phi = -\gamma^2 W (\gamma x_1 , x_2,x_3 )- \gamma^4 x_1^2 \underbrace{ \im \frac{ V\sin \alpha }{\left( (x_2 + ix_3)^2 - \gamma^2 x_1^2 \right)^{\frac32} } }_{=: Y ( \gamma x_1 , x_2,x_3)}
\]
Note that
\begin{align}
\| W (r, \cdot ) \|_{L^p (\R^2 )} &= C_p r^{-1+\frac2p } V \sin \alpha , \label{W_Lp_ests} \\
\| Y (r, \cdot ) \|_{L^q (\R^2 )} &= C_q r^{-3+\frac2q } V \sin \alpha \label{V_Lp_ests}
\end{align}
for all $p\in (2,4)$, $q\in [1,4/3)$, $r>0$. 

However, $\p_1 \phi$ has a jump at $x_1 = c$, so need to keep in mind that $\p_{11} \phi $ is a distribution,
\eqnb\label{distr_form}
\left\langle \p_{11} \phi , f \right\rangle  = \int_0^c \int_{\R^2}  \p_{11} \phi  f \, \d (x_2,x_3) \,\d x_1 -  \int_{\R^2} \p_1 \phi (\gamma c^-  , x_2,x_3) f(c,x_2,x_3) \,\d (x_2,x_3)
\eqne
for any $f\in C_c^\infty ( \R^3 )$. 

We now neglect issue \eqref{issue1}, that is pretend that \eqref{W_Lp_ests} holds also for $p=2$ and \eqref{V_Lp_ests} holds for $q=4/3$. Then  $\p_{11} \phi \in H^{-1/2} $, since
\[\begin{split}
| \left\langle \p_{11} \phi , f \right\rangle  | &\lec \gamma^2  c^{\frac12} V\sin \alpha  \| f  \|_{L^2 }+  \gamma^4 \int_0^c x_1^2 \| Y (\gamma x_1 , \cdot ) \|_{L^{4/3} (\R^2 )} \, \| f(x_1,\cdot )\|_{L^4 (\R^2 )} \d x_1\\
& \hspace{2cm} + \gamma^2 c \| W (\gamma c , \cdot )\|_{L^2 (\R^2 )} \| f (c,\cdot ) \|_{L^2 (\R^2 )}  \\
& \lec V \sin \alpha \left( \gamma^{\frac74} \| f \|_{L^2}+   \gamma^{\frac52 }\left( \int_0^c x_1 \,\d x_1 \right)^{\frac12} \| f \|_{L^2_{x_1} L^{4}_{x_2,x_3} } + \gamma^{\frac32} \| f (c,\cdot ) \|_{L^2 (\R^2 )}   \right) \\
&\lec V \sin \alpha \, \gamma^{\frac32} \| f \|_{H^{\frac12}}
\end{split} 
\]
for every $f\in C_c^\infty (\R^3)$, where in the last step we also neglected the borderline trace operator issue \eqref{issue2}, that is we (falsely) assumed that  $\| f (c,\cdot ) \|_{L^2 (\R^2 )}  \lec \| f \|_{H^{\frac12}}$. Hence
\[
\| \p_{11} \phi \|_{H^{-1/2}} \lec V \sin\alpha \, \gamma^{\frac32},
\]
and so \eqref{the_elliptic_issue} (which we take for granted here) gives that  $ \| \phi - \Phi \|_{H^{3/2} (\R^3 \setminus (\es \cup \widetilde{\es} ))} \lec V \sin\alpha \, \gamma^{\frac32} $, as required.

\subsection{Proof of Corollary~\ref{C1}}\label{sec_pf_cor}

Here we prove Corollary, that is we justify the expansions \eqref{Kp_and_Kv_claim}. Recalling (from~\eqref{gamma_def}) that $A=4\gamma $
\eqnb\label{temp11}
\begin{split}
|L_p^{(2D)} - L_p^{(3D)}| &= \rho V \cos^2 \alpha \left| \int_{\es} \p_1 (\phi - \Phi ) \right|\\
& \text{``}\lec\text{''} \rho V \cos^2 \alpha  \| \nabla (\phi - \Phi ) \|_{H^{1/2 }} \text{``}\lec\text{''} \rho V^2 \sin \alpha \cos^2 \alpha \, A^{\frac32},
\end{split}
\eqne
where, as mentioned in the introduction, we neglected issues \eqref{issue1}--\eqref{issue2} (which we pointed out using the quotation marks). Thus, recalling the definition~\eqref{Kp_def} of $K_p$ and recalling~\eqref{Kp_2d} that $K_p^{(2D)}=\frac{\pi}2A$,
\[
\left|K_p - \frac{\pi }2 A \right| = \frac{\left| C_{L_p}^{(3D)} - C_{L_p }^{(2D)}\right| }{\sin \alpha \cos^2 \alpha }  =  \frac{\left| L_p^{(2D)} - L_p^{(3D)}\right| }{\frac12 \rho V^2 |\es | \sin \alpha \cos^2 \alpha }  \lec \frac{1}{|\es |} A^{\frac32},
\]
which shows the first asymptotic expansion in~\eqref{Kp_and_Kv_claim}. The second one can be obtained by using Proposition~\ref{P1} to justify the modelling assumptions discussed in Section~\ref{sec_polhamus_explanation}. Namely, \eqref{Kv_delta} and the fact that $\cos \Lambda = \frac{A}4 \cos \beta = \frac{A}4 (1+O(A^2)) $ give
\[
\begin{split}
K_v &=  \left( K_p - \frac{K_p^2 }{\pi A} \right) \frac{1}{\cos \Lambda } = \left( \frac{\pi}2 A + O(A^{\frac32} ) - \frac{1 }{\pi A}\left( \frac{\pi}2 A + O(A^{\frac32} ) \right)^2  \right) \frac{4}{A (1+ O(A^2)) } \\
& =  \left( \frac{\pi}2  + O(A^{\frac12} ) - \frac{1 }{\pi A^2 }\left( \frac{\pi^2}4  A^2  + O(A^{\frac52} ) \right)  \right)  (4+ O(A^2)) \\
& =  \left( \frac{\pi}4  + O(A^{\frac12} \right)  (4+ O(A^2)) \\
&= \pi  + O(A^{\frac12})
\end{split}
\]
as $A\to 0$, justifying the second asymptotic expansion in \eqref{Kp_and_Kv_claim}.

\begin{rem}[about the borderline issues]\label{rem_borderline}
We conjecture that issues \eqref{issue1}--\eqref{issue2} (which were neglected above) can resolved using a more well-adapted functional setting. However, it is unlikely that the error $\phi - \Phi$ can be controlled in a topology significantly stronger than $H^{3/2}$. For example, such as $H^2$, as $\p_{11} \phi \not \in L^2 $. Perhaps a better  approximation  than $\phi$ could be used.

Moreover the appearance of the borderline nonintegrability issues in \eqref{temp11} cannot be resolved by any easy tricks. For example, replacing $H^s$ by $W^{s,p}$ does not work. Indeed, suppose that we want to estimate $\| \phi - \Phi \|_{W^{s,p}}$ instead of $H^{3/2}$. Then, in order to control the lift error (that is the first inequality in the second line of \eqref{temp11}) we need that $s>1+1/p$, as then $|\int_\es (\p_1 \phi - \p_1 \Phi )| \lec \| \phi - \Phi \|_{W^{s,p}}$. By the elliptic system $\Delta (\phi - \Phi ) = \p_{11} \phi$, this means that we need to estimate $\p_{11} \phi \in W^{s-2,p}$. However, estimating the second term on the right-hand side of \eqref{distr_form}, we get 
\[
\begin{split}
\left| {}_{W^{s-2,p}}\left\langle  \p_{11} \phi , f \right\rangle_{W^{2-s,p'}} \right| &\leq  \left| \int_0^c \int_{\R^2} \p_{11} \phi \, f \right |  + \gamma^2 c \int_{\R^2 } W(\gamma c , \cdot ) f(c,\cdot ) \\
&\leq \left| \int_0^c \int_{\R^2} \p_{11} \phi \, f \right |  + C_p \gamma^2 c  r^{-1+\frac2p } V \sin \alpha \, \| f (c, \cdot ) \|_{L^{p'}}
\end{split}
\]
for $p\in (2,4)$, where $p'\coloneqq p/(p-1)$. We can get that $\p_{11} \phi \in W^{s-2,p}$ if  the last norm of $f$ can be bounded by $\| f \|_{W^{2-s,p'}}$. However, this is true  only if $2-s >1/p'$, that is $s<1+1/p$, which is opposite to the above inequality. Thus the issue is analogous to the case when $p=2$.

\end{rem}

\section*{Acknowledgements}
This work was supported by the NSF grant no. DMS-2511556 and the Simons grant SFI-MPS-TSM-00014233. The author is grateful to  In-Jee Jeong, Anna Mazzucato, Daniel Peralta-Salas for interesting discussions.

\appendix

\section{2D point vortex suction}\label{app_2d_vortex_suc}

\begin{figure}[htbp]
\centering
    \includegraphics[width=9cm]{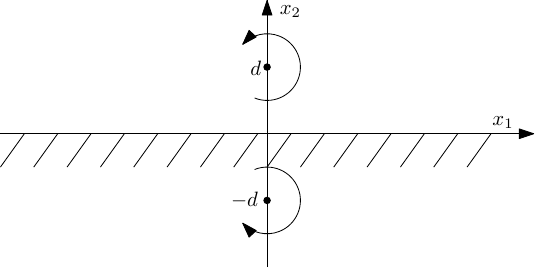}
  \caption{A sketch of the point vortex suction.}   \label{fig_2d_sucks}
\end{figure}

Here we follow \cite[Section~2.3.14]{nazarenko} to compute the suction force of a point vortex on a boundary. Consider a point vortex of strength $\kappa$, located at $(0,d)$, in the presence of a slip-boundary at $\{ (x_1,0) \colon x_1 \in \R \}$. In order to satisfy this boundary condition we consider an artificial point vortex at $(0,-d)$ of strength $-\kappa$, see Fig.~\ref{fig_2d_sucks}.

By the Biot-Savart law, the velocity field becomes
\[
\begin{split}
v (x) &= \frac{\kappa }{2\pi  (x_1^2 + (x_2-d )^2 )} \left( - (x_2-d) , x_1 \right) - \frac{\kappa }{2\pi  (x_1^2 + (x_2+d )^2 )} \left( - (x_2+d) , x_1 \right)   
\end{split}
\]
for $(x_1,x_2)\in \R^2_+$. For $x_2=0$ this becomes
\[
v(x_1,0^+) = \left(\frac{\kappa }{\pi (x_1^2 + d^2)}, 0 \right) .
\]
By the  steady Bernoulli equation $\frac{|v|^2 }2 + p = \mathrm{const}$ we can compute the pressure induced by the point vortex on the top of the $x_1$ axis,
\[
p_{\rm top} (x_1) = p_0 - \frac{\kappa^2 d^2}{\pi^2 (x_1^2 + d^2)^2},
\]
where $p_0$ denotes the limit pressure as $x_1\to \infty$. Thus the suction force that the point vortex exerts on the $x_1$ axis (in the direction of $x_2$) is
\[
\int_{-\infty}^\infty (p_0 - p_{\rm top}) \d x_1 = \frac{\kappa^2d^2 }{\pi^2} \int_{-\infty }^\infty \frac{\d x_1}{(x_1^2 + d^2)^2}  =  \frac{\kappa^2 }{2 \pi d }.
\]
We note in passing that an analogous computation as above can be used to explain how the weight of a flying airplane is transferred onto the ground, see~\cite[Section~17]{prandtl_NACA} for details.

\section{Potential flow past a line}\label{app_2dpot}
Here we briefly show that 
$-i  \sqrt{ (x+iy)^2-\frac{b^2}{4}  }$ is the complex potential of a potential flow past a flat horizontal line of length $b$,  which verifies the comment below \eqref{phi_def}.

Let $R\coloneqq \C^2 \setminus \{ |x|\leq b/2, y=0 \}$.
We recall the Joukovski map
\[
F(z) \coloneqq \frac{b}4 \left( z+\frac1z \right),
\]
which is a ``rotationless'' conformal mapping of the exterior of the unit disc, $\C \setminus \overline{D(0,1)}$ onto $R$. Thus, since $g(z) \coloneqq z$ is the complex potential of the uniform shear flow $(1,0)=\nabla \re g$, we have that $(g\circ F )(z) = F(z)$ is the complex potential of a potential flow past a unit disc in the direction of the positive real axis (since the composition of conformal mappings is a conformal mapping). In order to change the direction of the flow, we can multiply $z$ by $\ee^{i\theta}$ for a desired $\theta \in (0,2\pi)$. For example, $F(-iz )$ is the complex potential of a potential flow past a unit disc in the direction of the positive imaginary axis. Thus, since $F^{-1}$ maps $R$ onto $\C \setminus \overline{D(0,1)}$, we have that 
\[
F\left( -i F^{-1} (z) \right)
\]
is the complex potential of a flow in $R$ in the direction of the positive imaginary axis. A direct computation gives that\footnote{ We note that we wrote $\sqrt{\frac2b  z-1}\sqrt{\frac2b z+1}$ instead of $\sqrt{\frac{4}{b^2}z^2 -1 }$, since the latter results in an inevitable conflict regarding the branch choice. Namely, if we consider the branch cut of $z\mapsto \sqrt{z}$ along the negative real axis, then there exists a choice of $z\in R$ for which $F \left(\frac{2}{b} z + \sqrt{\frac{4}{b^2}z^2 -1 } \right) \ne z$. A similar problem appears when taking the positive real axis as the branch cut. Separating the two factors in the square root resolves this problem in both choices of the branch cut. This is fundamentally due to the fact that adding the real constant $\pm 1$ will never result in crossing the branch cut. }
\eqnb
F^{-1} (z) = \frac{2}bz + \sqrt{\frac2b  z-1}\sqrt{\frac2b z+1}
\eqne
and consequently
\eqnb
F\left( -i F^{-1} (z)\right) = - i \sqrt{z-\frac{b}2} \sqrt{z+\frac{b}2}.
\eqne

\section{2D airfoil of the flat line}

The circulation is 
\eqnb\label{circ_2D_flat_line}
\Gamma = \pi c V \sin \alpha 
\eqne

\section{The standard Munk lifting-line model}\label{app_munk}

Here we introduce some notions of the standard Munk lifting-lin model of $3$D airfoils. Let us represent a 3D wing by a continuous vortex line that is bound to the wing. To be precise, we consider a line
\[
\left\lbrace (0,y, 0) \colon |y |\leq \frac{b}2 \right\rbrace,
\]
where $b > 0$ denotes the wingspan, and we describe the aerodynamic property of the wing by defining the circulation $\Gamma (y)$ (of units $[m^2/s]$) as each spanwise location $y\in [-b/2,b/2]$. By the Kutta-Joukovski Theorem, if the lifting-line is placed in the background velocity field $(V,0,0)$, then the lift per unit span at a given spanwise location $y$ is $\rho V\, \Gamma (y)$, where $\rho$ denotes the density of the ambient fluid, so that the total lift is
\eqnb
L = \rho V \int_{-b/2}^{b/2} \Gamma (y) \, \d y.
\eqne

In terms of the lift coefficient $C_L \coloneqq \frac{L}{\frac12 \rho V^2 {|\es|}}$ this becomes
\eqnb
C_L = \frac{2}{V{|\es|}} \int_{-\frac{b}2}^{\frac{b}2} \Gamma = \frac{2 A}{Vb^2} \int_{-\frac{b}2}^{\frac{b}2} \Gamma ,
\eqne
where we substituted the aspect ratio
\eqnb\label{aspect_ratio}
A \coloneqq \frac{b^2}{|\es|}
\eqne
in the second equality.\\

We define the \emph{total effective circulation} $\Gamma_{\rm eff} $ as the average of $\Gamma(y)$,
\eqnb\label{Gammaeff}
\Gamma_{\rm eff} \coloneqq \frac1{b} \int_{-\frac{b}2}^{\frac{b}2} \Gamma (y) \, \d y.
\eqne
Thus, for a fixed lifting-line mode, $\Gamma_{\rm eff} $ represents the circulation that would need to be prescribed at each spanwise location of a rectangular wing (that is a lifting line with constant ciruclation) to obtain the same total lift.\\

Given $\Gamma (y)$ the bound vortex line gives rise to the trailing vortex sheet, by the Divergence Theorem. Indeed, since $\Gamma (y)$ varies in $y$, the divergence-free condition enforces that the vorticity must be shed, the modelling assumption is that this happens via an $e_1$ component of the vorticity vector, concentrated on the trailing sheet.

Using the 3D Biot-Savart law, $u [ \omega ] (x) =\frac{1}{4\pi } \int_{\R^3} \frac{(x-y)\times \omega (y) }{|x-y|^3} \d y $, we can compute the velocity field generated by the trailing vortex sheet. Restricting the velocity field to the lifting line, we have that
\eqnb\label{downwash_formula}
w[\Gamma ] (y) = \frac1{4\pi } \mathrm{p.v.} \int_{-\frac{b}2}^{\frac{b}2} \Gamma' (y') \frac{\d y' }{y-y'}
\eqne
(see~\cite[(5.15)]{anderson}) and then the induced drag can be computed using the Kutta-Joukovski theorem applied with velocity $w [\Gamma ]$,
\eqnb\label{ind_drag}
D_i [\Gamma ] = \rho \int_{-\frac{b}2}^{\frac{b}2} \Gamma (y) w[\Gamma ] (y) \d y.
\eqne

At this point it becomes advantageous to introduce the Chebyshev polynomials of the second kind, as observed by the author in \cite{ozanski_amo}. Namely, we write $\Gamma$ 
\eqnb\label{Gamma_exp}
\Gamma (y) = \sqrt{1-\left( \frac{2}b y \right)^2} \sum_{ k\geq 0} a_{2k} U_{2k} \left( \frac{2}b y \right) ,
\eqne
where $U_n(y) $ denotes the $n$-th Chebyshev polynomial of the second kind\footnote{recall~\cite{abramowitz_stegun} that $U_0 (x) \coloneqq 1$, $U_1 (x) \coloneqq 2x$ and $U_{n+1} (x) \coloneqq 2x U_n (x) - U_{n-1} (x)$ for $n\geq 1$.}. Note that we only consider even-numbered Chebyshev polynomials, since we assume that $\Gamma$ is even. Remarkably, the $U_{2k}$'s are eigenfunctions of the downwash operator $w[\Gamma ]$ in the sense that
\eqnb\label{cheb_efcn_fact}
w \left[ \sqrt{1-\left( \frac{2}b y \right)^2}   U_{2k}\left( \frac{2}b y \right) \right]  (x) = \frac{2k+1}{2b } U_{2k} \left( \frac{2}b x \right)
\eqne
for all $k\geq 0$, $x\in (-b/2,b/2)$, which follows immediately  from some well-known properties of Chebyshev polynomials\footnote{Some more classical works (see \cite{philips_hunsaker_56b, philips_hunsaker_joo_56a, philips_hunsaker_taylor} and references therein) make use of this fact in terms of the change of variable $\theta \coloneqq \cos (-2 x_2/b )$  and trigonometric identities.}, see \cite[(15)]{ozanski_amo} for details. Thus, for $\Gamma$ of the form \eqref{Gamma_exp},
\eqnb\label{wGamma_exp}
w [\Gamma ] (y) = \sum_{k\geq 0 } a_{2k} \frac{2k+1}{2b} U_{2k} \left( \frac{2}b y \right).
\eqne
Moreover, due to the fact that
$\{ U_{n} \}_{n\geq 0} $ form an orthonormal basis of the  $L^2 $ space on $(-1,1)$ with weight $\sqrt{1-(\cdot )^2}$, with  
\eqnb\label{Un_ortho}
\int_{-1}^1 U_n (\xi ) U_m (\xi ) \sqrt{1-\xi^2} \d \xi = \frac{\pi }2 \delta_{mn},
\eqne
we can write the induced drag \eqref{ind_drag} as 
\eqnb\label{ind_drag1}
\begin{split}
D_i[\Gamma ] &= \rho \sum_{m,k\geq 0}\frac{2m+1}{2b} a_{2k}a_{2m}  \int_{-\frac{b}2}^{\frac{b}2} U_{2k} \left( \frac{2}b y \right)  U_{2m }\left( \frac{2}b y \right)    \sqrt{1- \left( \frac{2}b y \right)^2 } \d y\\
& = \frac{\rho \pi}2 \sum_{k\geq 0}  \frac{2k+1}4a_{2k}^2,
\end{split}
\eqne
so that the induced drag coefficient
\[
C_{D_i} \coloneqq \frac{ D_i[\Gamma ] }{\frac12 \rho V^2 {|\es|} } = \frac{ \pi A }{ V^2 b^2} \sum_{k\geq 0}  \frac{2k+1}4a_{2k}^2,
\]
where we also recalled the definition \eqref{aspect_ratio} of the aspect ratio $A$.

On the other hand, the lift coefficient $C_L$ becomes
\[
C_L = \frac{2A}{Vb^2} \sum_{k\geq 0} a_{2k} \int_{-\frac{b}2}^{\frac{b}2} \sqrt{1-\left( \frac{2}b y \right)^2}  U_{2k} \left( \frac{2}b y \right) \d y = \frac{A\pi }{2 Vb } a_0,
\]
due to \eqref{Un_ortho}. Defining 
\[
K_i \coloneqq C_{D_i}/  C_L^2 
\]
we thus have 
\eqnb\label{Ki_comp}
K_i = \frac{\frac{ \pi A }{ V^2 b^2} \sum_{k\geq 0}  \frac{2k+1}4a_{2k}^2}{\frac{A^2\pi^2 }{4 V^2b^2 }a_0^2 } = \frac{1}{A\pi e}, 
\eqne
where
\[
e\coloneqq \frac{a_0^2}{\sum_{k\geq 0 } (2k+1) a_{2k}^2}
\]
denotes the efficiency coefficient. Note that $e\leq 1$, and $e=1$ if and only if $a_{2k}=0$ for all $k\geq 1$. In other words, we obtain the smallest induced drag to lift ratio (in the sense of minimizing $K_i$) only for the elliptic circulation distribution, 
\eqnb\label{Gamma_optimal}
\Gamma (y) = a_0 \sqrt{1-\left( \frac{2}b y \right)^2 }.
\eqne
Moreover, then
\eqnb\label{wGamma_optimal}
w[\Gamma ] (y) \quad \text{ is constant  and equal to }\,\, \frac{a_0}{2b} ,
\eqne
due to \eqref{wGamma_exp} and \eqref{wGamma_exp}. These last facts were first observed by Prandtl in his 1918 paper \cite{prandtl_trag}, who only considered the first two modes in the expansion \eqref{Gamma_exp}, and, instead of using Chebyshev polynomials, he derived a version of \eqref{wGamma_exp} in terms of polynomial basis $\{ 1 , y^2, y^4 , \ldots \}$  by employing the identities\footnote{These two identities replace the eigenfunction fact~\eqref{cheb_efcn_fact} for Chebyshev polynomials $U_{2k}$'s and the Chebyshev orthogonality property~\eqref{Un_ortho}, respectively.}
\[
\int_{-1}^1 \frac{\eta^{2n+1} \d \eta }{\sqrt{1-\eta^2} (\xi - \eta )} = \pi \sum_{m=0}^n p_{n-m} \,\xi^{2m}
\]
and
\[
\int_{-1}^1 \eta^{2n} \sqrt{1-\eta^2} \d \eta = \pi \,q_n
\]
for all $n\geq 0$, $\xi \in (-1,1)$
(which Prandtl~\cite{prandtl_trag} attributes to Albert Betz~\cite{betz_17}), where $p_n \coloneqq (1\cdot 3 \cdot 5 \cdot \ldots \cdot (2n-1) )/(2\cdot 4 \cdot \ldots \cdot 2n )$, and $q_n \coloneqq p_n / (2n+2)$.

\begin{rem}[The relevance of the assumption of fixed  $b$]
We emphasize that the above minimization result is  true only under the assumption that the wingspan $b>0$ is fixed (which was our assumption throughout this section). In the case when $b$ becomes a variable, the minimization problem becomes nonconvex (see the comments in~\cite[p.~39]{ozanski_amo}), and so becomes much more challenging.
\end{rem}
\begin{rem}[More efficient solutions if $b>0$ is a variable]
 In the case when $b$ is assumed to be a variable of the problem of minimizing induced drag, there are two natural modelling assumptions: (1) that the wing has a fixed weight, and (2) that the area of the wing's intersection at any fixed spanwise location must be big enough to support the bending moment resulting from the lift force acting at points from the location until the wingtip. Such modelling assumptions were used in the celebrated 1933 paper of Prandtl~\cite{prandtl_33} (see also a recent translation~\cite{prandtl_33_trans} by Hunsaker and Philips), who showed that the optimal solution in such case is a bell-shaped circulation distribution, which is $11\%$ more efficient. The latter paper provides a clear explanation why birds' wings do not have elliptic planforms, and was explored experimentally by NASA's Armstrong Flight Research Center~\cite{bmjeg,newton}, see also the series of works on wing design by Hunsaker, Philips, and collaborators~\cite{philips_hunsaker_56b, philips_hunsaker_joo_56a, philips_hunsaker_taylor}. We also note a recent improvement \cite{ozanski_amo} of Prandtl's 1933 result~\cite{prandtl_33}, which, except for the bending moment due to the lift force, also takes into account the bending moment due to gravity, which provides an even more efficient solution.
\end{rem}

\bibliographystyle{plain}
\bibliography{literature}

\end{document}